\documentclass[12pt]{article}
\usepackage{amsfonts}
\usepackage{amssymb}
\usepackage{polski}

\input{tcilatex}
\begin{document}

\begin{center}
\textbf{Additivity, subadditivity and linearity:}

\textbf{automatic continuity and quantifier weakening\\[0pt]
by \\[0pt]
N. H. Bingham and A. J. Ostaszewski}\\[0pt]

\bigskip

\textit{To Charles Goldie.}

\bigskip
\end{center}

\noindent \textbf{Abstract. }We study the interplay between additivity (as
in the Cauchy functional equation), subadditivity and linearity. We obtain
automatic continuity results in which additive or subadditive functions,
under minimal regularity conditions, are continuous and so linear. We apply
our results in the context of quantifier weakening in the theory of regular
variation, completing our programme of reducing the number of hard proofs
there to zero.

\bigskip



\section{Introduction}

Our main theme here is the interplay between additivity (as in the Cauchy
functional equation), subadditivity and linearity. As is well known, in the
presence of smoothness conditions (such as continuity), additive functions $%
A:\mathbb{R}\rightarrow \mathbb{R}$ are linear, so of the form $A(x)=cx$.
There is much scope for \textit{weakening} the smoothness requirement and
also much scope for weakening the universal quantifier by\textit{\ thinning}
its range $\mathbb{A}$ below from the classical context $\mathbb{A}=\mathbb{R%
}$:%
\begin{equation}
A(u+v)=A(u)+A(v)\text{ }(\forall u,v\in \mathbb{A}). 
\tag{\QTR{rm}{Add}$_{\QTR{Bbb}{A}}$($A$)}
\end{equation}%
We address the Cauchy functional equation in \S 2 below. The philosophy
behind our quantifier weakening\footnote{%
or `quantifier easing'} theorems is to establish \textit{linearity} of a
function $F$ on $\mathbb{R}$ from its \textit{additivity} on a thinner set $%
\mathbb{A}$ and from additional (`side') conditions, which include its 
\textit{extendability} to a subadditive function; recall that $S:\mathbb{R}%
\rightarrow \mathbb{R}\cup \mathbb{\{-\infty },\mathbb{+\infty \}}$ is 
\textit{subadditive }[HilP, Ch. 3] if for all $u,v\in \mathbb{R}$%
\begin{equation}
S(u+v)\leq S(u)+S(v),  \tag{\QTR{rm}{Sub}}
\end{equation}%
whenever meaningful on the right-hand side (cf. [Kuc, 16.1] for subbadditive
and [Roc, p. 23 ffol] for convex functions); see also [MatS], [BinO4]. (This
choice for the setting is more convenient than alternatively working on $%
\mathbb{R}_{+}$, even though one-sided side-conditions are important here.)
To motivate our main result we begin with an automatic continuity theorem,
devoted entirely to subadditive functions; it implies a result about those
linear functions that have subadditive extensions -- see Proposition 7 below
(on uniqueness of extension). This (Theorem 0 below) makes explicit an
argument springing from a step in a proof by Goldie, of Th. 3.2.5 in [BinGT]
(BGT below, for brevity), recently improved and generalized in [BinO14]
(though still implicit even there).

We recall that for $S$ subadditive and finite-valued, $S(0)\geq 0,$ as $%
S(0)\leq S(0)+S(0),$ so that $S(0)=0$ iff $S(-z)=-S(z)$ for some $z,$ as
will be the case below when $S\ $extends an additive function; cf. [Kuc, p.
457].

\bigskip

\noindent \textbf{Theorem 0. }\textit{For subadditive} $S:\mathbb{R}%
\rightarrow \mathbb{R}\cup \{-\infty ,+\infty \}$ \textit{with }$%
S(0+)=S(0)=0:$\textit{\ }$S$\textit{\ is continuous at }$0$ \textit{iff }$%
S(z_{n})\rightarrow 0,$\textit{\ for some sequence }$z_{n}\uparrow 0,$%
\textit{\ and then }$S$ \textit{is continuous everywhere, if finite-valued.}

\bigskip

The last part above draws on [HilP, Th. 2.5.2] that, for a subadditive
function, continuity at the origin implies continuity everywhere. Theorem 0
above, in the presence of right-sided continuity, asserts that the merest
hint of left-sided continuity gives full continuity; contrast this with the
behaviour of the subadditive function $\mathbf{1}_{[0,\infty )}$, which is
continuous on the right but not on the left. This leads to the question of
whether right-sided continuity can be thinned out. We are able to do so in
the next two results below, but at the cost of imposing more structure,
either on the left, or on the right. We need the following two definitions.

\bigskip

\noindent \textbf{Definitions.} 1. Say that $\Sigma $ is \textit{locally
Steinhaus-Weil (SW)}, or has the \textit{SW property locally}, if for $%
x,y\in \Sigma $ and, for all $\delta >0$ sufficiently small, the sets%
\[
\Sigma _{z}^{\delta }:=\Sigma \cap B_{\delta }(z), 
\]%
for $z=x,y,$ have the \textit{interior-point property,} that $\Sigma
_{x}^{\delta }\pm \Sigma _{y}^{\delta }$ has $x\pm y$ in its interior. (Here 
$B_{\delta }(x)$ is the open ball about $x$ of radius $\delta .)$ See
[BinO9] for conditions under which this property is implied by the
interior-point property of the sets $\Sigma _{x}^{\delta }-\Sigma
_{x}^{\delta }$ (cf. [BarFN]); for a rich list of examples, see \S 4.

\noindent 2. Say that $\Sigma \subseteq \mathbb{R}$ is \textit{shift-compact 
}if for each \textit{null sequence} $\{z_{n}\}$ (i.e. with $z_{n}\rightarrow
0$) there are $t\in \Sigma $ and an infinite $\mathbb{M\subseteq N}$ such
that%
\[
\{t+z_{m}:m\in \mathbb{M}\}\subseteq \Sigma . 
\]%
See [BinO11], and for the group-action aspects, [MilO].

\bigskip

For connections between these two, and related result, see \S 8.3. Thus
armed, we begin with symmetric (two-sided) thinning. The second part of the
result below is a variant of the Berz Theorem for \textit{sublinear}
functions ([Berz], \S 5, [BinO16]).

\bigskip

\noindent \textbf{Theorem 0}$^{\prime }$\textbf{\ }[BinO14, Th. 3]\textbf{.} 
\textit{If }$S:\mathbb{R}\rightarrow \mathbb{R}$\textit{\ is subadditive
with }$S(0)=0$ \textit{and there is a symmetric set }$\Sigma $\textit{\
containing }$0$ \textit{with:}

\noindent (i)\textit{\ }$S$\ \textit{is continuous at }$0$\textit{\ on }$%
\Sigma $\textit{;}

\noindent (ii) \textit{for all small enough }$\delta >0,$ $\Sigma
_{0}^{\delta }$ \textit{is locally Steinhaus-Weil}\newline
\noindent -- \textit{then }$S$ \textit{is continuous at }$0$\textit{\ and so
everywhere. }

\textit{In particular, this conclusion holds if there is a symmetric set }$%
\Sigma $\textit{, Baire/measurable and non-negligible in each }$(0,\delta )$%
\textit{\ for }$\delta >0,$\textit{\ on which}%
\[
S(u)=c_{\pm }u\text{ }\mathit{for\ some\ }c_{\pm }\in \mathbb{R}\mathit{\
and\ all\ }u\in \mathbb{R}_{+}\cap \Sigma \mathit{,or}\text{ }\mathit{all}%
\text{ }u\in \mathbb{R}_{-}\cap \Sigma \text{ }\mathit{resp.} 
\]

\bigskip

The alternative case to two-sided thinning is one-sided thinning accompanied
by linear bounding. Although one tries to impose continuity conditions on a
thin set below, these cannot be `too thin' as the example of the indicator
function of the irrationals shows: $\mathbf{1}_{\mathbb{R}\backslash \mathbb{%
Q}}$ is subadditive and additive on $\mathbb{Q}$ (indeed $\mathbb{Q}$%
-homogeneous), but not continuous. More thinning is possible by involving
more structure: the ability to `span' (see below).

To motivate the accompanying side-condition, consider a subadditive $S$ that
is locally bounded, say by some $\varepsilon >0$ on $(0,a]$. For any $x>0,$
choose $n\in \mathbb{N}$ with $(n-1)<x/a\leq n;$ then%
\[
S(x)\leq (n-1)S(a)+S(x-(n-1)a)\leq xS(a)/a+\varepsilon . 
\]%
In particular, for $x\geq a,$ $S(x)\leq c_{a}x,$ for $c_{a}:=[\varepsilon
a+S(a)]/a,$ i.e. $S\ $is linearly bounded away and to the right of the
origin. Theorem 0$^{+}$ below derives global linear boundedness from similar
behaviour on a thin set near $0$ when $c_{a}$ itself is bounded above on the
thin set.

\bigskip

\noindent \textbf{Theorem 0}$^{+}$\textbf{.} \textit{Let} $\Sigma \subseteq
\lbrack 0,\infty )$ \textit{be locally SW accumulating at }$0$\textit{.
Suppose }$S:\mathbb{R}\rightarrow \mathbb{R}$\textit{\ is subadditive with }$%
S(0)=0$ \textit{and:}

\noindent $S|\Sigma $ \textit{is linearly bounded above by }$G(x):=cx$%
\textit{, i.e.} $S(\sigma )\leq c\sigma $ \textit{for some }$c$\textit{\ and
all }$\sigma \in \Sigma ,$ \textit{so that in particular,}%
\[
\lim \sup\nolimits_{\sigma \downarrow 0,\text{ }\sigma \in \Sigma }S(\sigma
)\leq 0. 
\]%
\textit{\ Then }$S(x)\leq cx$ \textit{for all }$x>0,$\textit{\ so} 
\[
\lim \sup_{x\downarrow 0}S(x)\leq 0, 
\]%
\textit{and so }$S(0+)=0.$

\textit{In particular, if furthermore there exists a sequence }$%
\{z_{n}\}_{n\in N}$\textit{\ with }$z_{n}\uparrow 0$\textit{\ and }$%
S(z_{n})\rightarrow 0,$\textit{\ then }$S$\textit{\ is continuous at }$0$%
\textit{\ and so everywhere.}

\bigskip

Boundedness in the latter case could be provided on an open set $U\subseteq
(0,1)$ accumulating at $0;$ continuity in the former case on a set $\Sigma
=(-U)\cup \{0\}\cup U$ with $U$ as before. Indeed, in this context one may
equivalently assume that the set $\Sigma $ has precisely this form: see \S 4.

\bigskip

Theorem 1, our quantifier weakening theorem, was our original motivation,
for reasons explained later in the paper. The Heiberg-Seneta side-condition
of Theorem 1 is due to Heiberg ([Hei] in 1971) and Seneta ([Sen] in 1973) --
see BGT\ Th. 1.4.3.

\bigskip

\noindent \textbf{Theorem 1. }\textit{For} $S:\mathbb{R}\rightarrow \mathbb{R%
}\cup \mathbb{\{-\infty },\mathbb{+\infty \}}$ \textit{subadditive, and }$%
\mathbb{A}\subseteq \mathbb{R}$\textit{\ an additive subgroup, suppose that}%
\newline
\noindent (i) $\mathbb{A}$\textit{\ is dense;}\newline
\noindent (ii) $A:=S|\mathbb{A}$\textit{\ is finite and additive, i.e. }$%
\mathrm{Add}_{\mathbb{A}}$($S$) \textit{holds};\newline
\noindent (iii)\textit{\ }$S$\textit{\ satisfies the one-sided
(Heiberg-Seneta) boundedness condition}%
\begin{equation}
\lim \sup_{u\downarrow 0}S(u)\leq 0.  \tag{$HS(S)$}
\end{equation}%
\textit{\ Then }$S$\textit{\ is linear: }$S(u)=cu$\textit{\ for some} $c\in 
\mathbb{R}$\textit{, and all }$u\in \mathbb{R}.$

\bigskip

Its proof relies on Theorem 0. Our next result, Theorem 1$^{\prime }$ based
(in part) on Theorem 0$^{\prime }$, is formulated in the spirit of Theorem 1
so far as the Heiberg-Seneta-style condition is concerned. However, the
passage to the limit below is through a set $\Sigma $ which may be (very)%
\textit{\ thin} (`ghost-like'); the limit may be one-sided or two-sided,
depending both on $\Sigma $ and the ambient context.

\bigskip

\noindent \textbf{Theorem 1}$^{\prime }$ (cf. [BinO15, \S 6 Th. 5]). \textit{%
Theorem 1 above holds with condition }(iii)\textit{\ replaced by any one of
the following:}\newline
\noindent (iii-a)\textit{\ }$S$\textit{\ satisfies the Heiberg-Seneta
boundedness condition thinned out to a symmetric set }$\Sigma $\textit{\
that is locally SW, i.e.}%
\[
\lim \sup\nolimits_{u\rightarrow 0,\text{ }u\in \Sigma }S(u)\leq 0; 
\]%
\noindent (iii-b)\textit{\ }$S$\textit{\ is linearly bounded above on a
locally SW subset }$\Sigma \subseteq $\textit{\ }$\mathbb{R}_{+}=(0,\infty )$%
\textit{\ accumulating at }$0,$\textit{\ so that in particular}%
\[
\lim \sup\nolimits_{u\downarrow 0,\text{ }u\in \Sigma }S(u)\leq 0; 
\]%
\textit{\ }\noindent (iii-c)\textit{\ }$S$\textit{\ is bounded above on a
locally SW subset }$\Sigma \subseteq $\textit{\ }$\mathbb{A}_{+}$\textit{\
accumulating at }$0$\textit{, that is, the following }$\lim \sup $ \textit{%
is finite:}

\begin{equation}
\lim \sup\nolimits_{u\downarrow 0,\text{ }u\in \Sigma }S(u)<\infty ; 
\tag{$SW$-$HS(S)$}
\end{equation}%
\noindent (iii-d) $S$\textit{\ is bounded on a subset }$\Sigma \subseteq $%
\textit{\ }$\mathbb{A}$\ \textit{that is shift-compact (e.g. on a set that
is locally SW, and so on an open set) and\ }%
\[
\mathbb{A}=\mathbb{A}_{S}:=\{u:G(u):=\lim_{x\rightarrow \infty }[S(u+x)-S(x)]%
\text{ exists and is finite}\}. 
\]

\bigskip

Theorem 1$^{\prime }$(c) above encompasses as an immediate corollary
Ostrowski's Theorem [BinO11] and its classical generalizations. Below we use 
\textit{negligible} to mean meagre or null, according as (Baire) category or
(Lebesgue) measure is considered; we use \textit{non-negligible} to mean a
Baire/Lebesgue set that is not correspondingly negligible. Further
discussion here involves set-theoretic assumptions (typically, alternatives
to the alternatives to the Axiom of Choice, AC, see \S 8.5 for references).

\bigskip

\noindent \textbf{Theorem O }([Darb], [Ostr], [Meh]; cf. [BinO11])\textbf{. }%
\textit{If }$A:\mathbb{R}\rightarrow \mathbb{R}$ \textit{is additive and
bounded above on a non-negligible set, then }$A$ \textit{is linear.}

\bigskip

The measure-theoretic development here came earlier chronologically; it was
then noticed that the Baire (or category) case was closely analogous. The
two cases are developed in parallel in BGT. It has emerged recently that the
primary case is in fact the category case; see e.g. [BinO2,10], [Ost2].

The point of introducing side-conditions of Heiberg-Seneta type is that they
enable one to avoid assuming that our functions are \textit{either Baire or
measurable}. This is in stark contrast to all the other major results in the
theory of regular variation (for which see below), when such an assumption
is necessary -- e.g., to avoid the \textit{Hamel pathology} (see e.g. [HilP,
p. 238], BGT p. 5) of discontinuous additive functions (necessarily wildly
non-measurable/non-Baire). See e.g. the property of local boundedness in
Prop. 9 below (\S 4), and foundational results such as the Uniform
Convergence Theorem [BinO3]. We strengthen results of this type by imposing
the side-condition on \textit{as thin a set as possible} (`ghost-like
Heiberg-Seneta conditions'), albeit with some regularity of structure.

Our results here concern, as well as quantifier weakening and automatic
continuity, various results on additive, subadditive and sublinear
functions. Our original motivation was (quantifier weakening in the context
of) regular variation (\S 7: Karamata theory and its extensions; see e.g.
BGT). This specific motivation turns out to be valuable here: our viewpoint
on the area, informed by this, is complementary to (indeed, contrasts with)
that of the standard work in the area, Kuczma ([Kuc, Ch. 16]).

The rest of the paper is structured as follows. The backdrop of the Cauchy
equation is considered in \S 2; then we pass in \S 3 to developing the theme
of linearity first from additivity then from subadditivity, proving Theorems
0 and 1 and their variants, and also Theorem 2. We stop to clarify the
thinning aspect of the classical \textit{Heiberg-Seneta condition} in \S 4.
Using the theme of linearity from subadditivity we extend in \S 5 Berz's
Theorem on sublinearity [Berz] from the classical measurable case to the
Baire category case, more natural for us here, since refining the Euclidean
topology to the density topology (which converts measurable sets into 
\textit{Baire }sets, i.e. sets with the Baire property -- see e.g. [BinO10,
17]) yields a parallel new proof of Berz's Theorem. In \S 6 we discuss
thinning by \textit{spanning}: we want to \textit{weaken quantifiers} by 
\textit{thinning their range} as much as possible. But limits are imposed on
this: a set which is too thin will not be able to span. We are thinking here
of the reals $\mathbb{R}$ as a vector space over the rationals $\mathbb{Q}$,
Hamel bases etc. All this is motivated by regular variation, for which see 
\S 7. This makes good on the claim we have already made elsewhere (see
[BinO7,8,11,12]): this reduces the number of hard proofs in the theory of
regular variation to zero.

\section{Cauchy theory}

Theorem 1 (to be proved in \S 3, where further results of this kind are
established) is concerned with \textit{weakening} the quantifier in the
classical Cauchy functional equation, by thinning the range $\mathbb{A}$ of
the universal quantifier in (\textrm{Add}$_{\mathbb{A}}$) above from the
classical context $\mathbb{A}=\mathbb{R}$. The philosophy behind the
theorems is to establish \textit{linearity} of a function $F$ on $\mathbb{R}$
from its \textit{additivity} on a thinner set $\mathbb{A}$ and from
additional (`side') conditions, which include its extendability to a
subadditive function. The principal motivation for such an approach arises
in regular variation (\S 7) and rests on the following result, which
identifies the \textit{additive kernel} $G$ of the function $F^{\ast }$
below. The paper studies conditions under which $F^{\ast }$ coincides with
this kernel. The condition (i) below motivated further study in [BinO15, \S %
5 especially Prop. 6] (within a more subtle group structure, and under the
more demanding requirements of uniform convergence).

\bigskip

\noindent \textbf{Proposition 1 (Additive Kernel). }\textit{For }$F:\mathbb{R%
}\rightarrow \mathbb{R}$\textit{\ put}%
\[
\mathbb{A}_{F}:=\{u:G(u):=\lim_{x\rightarrow \infty }[F(u+x)-F(x)]\text{
exists and is finite}\}, 
\]%
\textit{\ and, for }$u\in \mathbb{R}$\textit{\ define}%
\[
F^{\ast }(u):=\lim \sup\nolimits_{x\rightarrow \infty }[F(u+x)-F(x)]. 
\]%
\textit{Then:}\newline
\noindent (i)\textit{\ }$\mathbb{A}_{F}$\textit{\ is an additive subgroup;}%
\newline
\noindent (ii) $G$ \textit{is an additive function on }$\mathbb{A}_{F}$%
\textit{;}\newline
\noindent (iii) $F^{\ast }:\mathbb{R}\rightarrow \mathbb{R}\cup \{+\infty
,-\infty \}$\textit{\ is a subadditive extension of }$G$\textit{;}\newline
\noindent (iv) $F^{\ast }$\textit{\ is finite-valued and additive iff }$%
\mathbb{A}_{F}=\mathbb{R}$ \textit{and }$F^{\ast }(u)=G(u)$ \textit{for all }%
$u.$

\bigskip

\noindent \textbf{Proof. }(i) $0\in \mathbb{A}_{F}$. Next, from%
\begin{equation}
F(u+v+x)-F(x)=[F(u+v+x)-F(v+x)]+[F(v+x)-F(x)],  \tag{*}
\end{equation}%
we see that $\mathbb{A}_{F}$ is a subsemigroup of $\mathbb{R}.$ In fact it
is a subgroup: for $u\in \mathbb{A}_{F}$ one has $-u\in \mathbb{A}_{F},$
because on writing $y=u+x$ one has 
\[
F(-u+y)-F(y)=-[F(u+x)-F(x)].\text{ }\square \text{ (i)} 
\]%
\noindent (iii) The identity ($^{\ast }$) also implies subadditivity of $%
F^{\ast }$ and the partial additivity result that 
\begin{equation}
F^{\ast }(u+v)=F^{\ast }(u)+F^{\ast }(v)\text{ }\forall u\in \mathbb{A}_{F}%
\text{ }\forall v\in \mathbb{R}.  \tag{$^{\ast \ast }$}
\end{equation}%
We note that for $u\in \mathbb{A}_{F},$ i.e. when $G(u)$ exists, then $%
F^{\ast }(u)=G(u),$ so proving part (iii) for $u\in \mathbb{A}_{F}$. $%
\square $ (iii)

\noindent (ii) With $v\in \mathbb{A}_{F},$ (**) above yields additivity of $%
G $. $\square $ (ii)

\noindent (iv) If $F^{\ast }$ is finite-valued and additive, then $F^{\ast
}(-u)=-F^{\ast }(u)$ for all $u.$ The substitution $y=u+x$ yields%
\begin{eqnarray*}
\lim \inf\nolimits_{x\rightarrow \infty }[F(u+x)-F(x)] &=&-\lim
\sup\nolimits_{x\rightarrow \infty }[F(x)-F(u+x)] \\
&=&-\lim \sup\nolimits_{y\rightarrow \infty }[F(-u+y)-F(y)] \\
&=&-F^{\ast }(-u)=F^{\ast }(u).
\end{eqnarray*}%
That is, $\lim \inf\nolimits_{x\rightarrow \infty }[F(u+x)-F(x)]=\lim
\sup\nolimits_{x\rightarrow \infty }[F(u+x)-F(x)],$ i.e. $F^{\ast
}(u)=\lim_{x\rightarrow \infty }[F(u+x)-F(x)]=G(u)$ for all $u$; so $\mathbb{%
A}_{F}=\mathbb{R}.$ The converse follows from (i) and (ii). $\square $ (iv). 
$\square $

\bigskip

\noindent \textbf{Remark. }With the Axiom of Choice \textrm{AC} replaced by
an axiom (such as the Axiom of Determinacy, \textrm{AD}) under which all
sets are Baire/measurable (a price that most mathematicians most of the time
will not be willing to pay!), Prop. 1 (helped by Prop. 6 below) is all we
need, and the remaining results below become unnecessary. See [BinO16,
Appendix 1] and \S 8.5; cf. [BinO15, \S 7 Th. 6].

\bigskip

Proposition 2 below leads from Theorem 1 to Theorem 1$^{\prime }$ in \S 3
below.

\bigskip

\noindent \textbf{Proposition 2.} \textit{For }$A:\mathbb{R}\rightarrow 
\mathbb{R}$\textit{\ additive, the following are equivalent:\newline
}\noindent (i) $A$\textit{\ is bounded above on a non-negligible
Baire/measurable set;\newline
}\noindent (ii)\textit{\ }$A$ \textit{is bounded above on an interval;%
\newline
}\noindent (iii)\textit{\ for some} $M\in \mathbb{R}$%
\begin{equation}
\lim \sup_{u\downarrow 0}A(u)\leq M;  \tag{$\lim \sup_{M}$}
\end{equation}%
\noindent (iv) $(\lim \sup )_{0}$ \textit{holds, i.e. }$(HS(A))$\textit{\
holds;}

\bigskip

\noindent \textbf{Proof}. (i)$\rightarrow $(ii): If $A$ is bounded above on
a non-negligible Baire/measurable set $L$, then it is bounded above on $L+L,$
which contains an interval by the \textit{Steinhaus-Piccard-Pettis
Sum-Theorem} ([BinO9, Th. E]; cf. [GroE], [BinO11] and the recent [BinO18]).

\noindent (ii)$\rightarrow $(iii): By additivity, $(\lim \sup )_{M}$ holds
for some $M\in \mathbb{R}.$

\noindent (iii)$\rightarrow $(iv): Without loss of generality (w.l.o.g.) $%
M\geq 0.$ For any $K>M,$ if $\sup \{A(u):0<u<\delta \}<K,$ then by
additivity $\sup \{A(u):0<u<\delta /2\}<K/2$, as $2A(u)=A(2u)<K,$ so $(\lim
\sup )_{M/2}$ holds, and so the least $M\geq 0$ for which the condition
holds is $M=0.$

\noindent (iv)$\rightarrow $(i): Clear. $\square $

\bigskip

\textbf{Remark.} Once Theorem 1 is established, we may apply Theorem 1 to $%
S=A$ with $\mathbb{A}=\mathbb{R}$ to deduce a further equivalent condition:

\noindent (v) $A$\textit{\ is linear on }$\mathbb{R}$\textit{: }$A(u)=cu$%
\textit{\ for some} $c\in \mathbb{R}$\textit{, and all }$u\in \mathbb{R}.$

\bigskip

The analogue of Prop. 2 for an additive subgroup is also relevant below.

\bigskip

\noindent \textbf{Proposition 2}$^{\prime }$\textbf{\ (Automatic continuity,
after Darboux, }[Darb]\textbf{). }\textit{For }$\mathbb{A}$\textit{\ a dense
additive subgroup and }$A:\mathbb{A}\rightarrow \mathbb{R}$\textit{\
additive, the following are equivalent:}

\noindent (i) $A$\textit{\ is continuous};

\noindent (ii) $A$\textit{\ is right-continuous at} $0;$

\noindent (iii) $A$ \textit{is continuous at} $0;$

\noindent (iv) $A$ \textit{is locally bounded;}

\noindent (v) $A$ \textit{is locally bounded above on some interval.}

\bigskip

\noindent \textbf{Proof.} As this is routine, we refer to [BinO11] for
details, save to say that (v)$\rightarrow $(iv) is as in Prop. 5 below, and
to show that (iv)$\rightarrow $(i). If $A$ is locally bounded at $0,$ there
is $\delta >0$ and $M$ such that $|A(a)|\leq M$ for $a\in \mathbb{A}$ with $%
|a|<\delta .$ Given $\varepsilon >0,$ choose an integer $N$ with $%
N>M/\varepsilon .$ For $a\in \mathbb{A}$ with $|a|<\delta /N,$ $%
|NA(a)|=|A(Na)|\leq M,$ so%
\[
|A(a)|\leq \varepsilon , 
\]%
as $M/N<\varepsilon .$ This gives continuity.$\square $

\bigskip

\textbf{Remark.} In the next section Prop. 6 establishes the further
equivalent condition:

\noindent (vi) $A$\textit{\ is linear: }$A(u)=cu$\textit{\ for some} $c\in 
\mathbb{R}$\textit{, and all }$u\in \mathbb{A}.$

\section{Linearity from subadditivity}

Here we establish a paradigm for identifying circumstances \textit{when
linearity may be deduced from subadditivity} -- encapsulated in Theorem 2
below -- by showing that $S(y)/S(x)=y/x$ on a dense subspace and appealing
to right-continuity.

We begin with a result linking the Heiberg-Seneta condition $(HS)$ of
Theorem 1 with automatic one-sided continuity. We need some preliminaries:
we turn first to conditions implying finite-valued subadditivity. The first
result requires a mixture of one-sided and two-sided information.

\bigskip

\noindent \textbf{Proposition 3. }\textit{For} $S:\mathbb{R}\rightarrow 
\mathbb{R\cup \{-\infty },\mathbb{+\infty \}}$\textit{\ subadditive, write }$%
\Sigma _{+}:=\{u\in \mathbb{R}:S(u)<\infty \}.$ \textit{If }$\Sigma _{+}\cap
\lbrack 0,\infty )$ \textit{contains an interval and }$\Sigma
_{+}\nsubseteqq \lbrack 0,\infty )$\textit{, then either }$S$ \textit{is
finite everywhere or is identically }$-\infty .$ \textit{In particular, if}

\noindent (i)\textit{\ }$S$ \textit{is finite on a subset }$\Sigma $ \textit{%
unbounded below} \textit{(e.g. a dense subset of }$\mathbb{R}$\textit{);}

\noindent (ii)\textit{\ }$S$\textit{\ is bounded above on }$(0,\delta )$ 
\textit{for some }$\delta >0,$ \textit{e.g.} $S$ \textit{satisfies the
condition }$(HS(S)),$

\noindent -- \textit{then }$S$\textit{\ is finite everywhere.}

\bigskip

\noindent \textbf{Proof.} As $\Sigma _{+}$ is a subsemigroup of the additive
group $\mathbb{R}$, and contains an interval, it contains a ray $[A,\infty )$
(see e.g. BGT\ Cor. 1.1.5). Choose $c\in (-\infty ,0)\cap \Sigma _{+}.$ Then 
$nc+(A,\infty )\subseteq \Sigma _{+}$ for all $n\in \mathbb{N}$, so $\Sigma
_{+}=\mathbb{R}$, i.e. $S(u)<\infty $ for all $u\in \mathbb{R}.$ If $%
S\not\equiv -\infty ,$ say $S(u_{0})>-\infty ,$ then for any $u\in \mathbb{R}
$%
\[
S(u)\geq S(u_{0})-S(u_{0}-u)>-\infty . 
\]

\noindent The particular case now follows, since by (i) $\Sigma \nsubseteqq
\lbrack 0,\infty ),$ and by (ii) $S$ is bounded on an interval: there is $%
\delta >0$ such that $S(x)<1$ for all $x\in (0,\delta )$, so $\Sigma \cap
\lbrack 0,\infty )$ contains $(0,\delta ).$ $\square $

\bigskip

We will soon prove in Proposition 3$^{\prime }$ a one-sided variant. The
argument used can yield more; so, it is more convenient to prove first

\bigskip

\noindent \textbf{Proposition 4.} \textit{If} $S:\mathbb{R}\rightarrow 
\mathbb{R}\ $\textit{with }$S(0)=0$\textit{\ is subadditive and linearly
bounded above by }$G(x)=cx$\textit{\ on an open set }$U$\textit{\
accumulating to the right at }$0,$\textit{\ then }$S$\textit{\ is linearly
bounded above by }$G$\textit{\ on }$\mathbb{R}_{+}.$

\textit{Furthermore, if }$S(x)=G(x)$\textit{\ on a dense set }$D,$\textit{\
then }$S(x)=cx$\textit{\ on }$\mathbb{R}_{+}.$

\bigskip

\noindent \textbf{Proof.} The set $\Sigma _{+}:=\{v:S(v)\leq cv\}$ is an
(additive) semigroup containing $U$. By a theorem of Kingman [BinO18, Th.
3.5], $\Sigma _{+}$ is dense in $\mathbb{R}_{+}:$ for any interval $%
I\subseteq \mathbb{R}_{+}$ there is $\eta \in I$ such that $\eta /m\in U$
for infinitely many $m\in \mathbb{N};$ for such an $\eta $ and any
corresponding $m\in \mathbb{N}$,%
\[
S(\eta )\leq mS(\eta /m)\leq mc(\eta /m)=c\eta . 
\]%
So $\eta \in \Sigma _{+}\cap I,$ proving density. Since $\Sigma _{+}$ is a
semigroup $\eta +(a,b)\subseteq \Sigma _{+}$ for any interval $%
(a,b)\subseteq U,$ with $a$ chosen as small as desired, since $U\ $%
accumulates at $0.$ So the family of open intervals $J$ contained in $\Sigma
_{+}$ have dense union in $\Sigma _{+}$: in sum, \textrm{int}$(\Sigma _{+})$
accumulates to the right and left of any point of $(0,1).$ Fix any $x\in
(0,1).$ There exists $h>0$ as small as desired such that $x-h\in \Sigma
_{+}. $ For such an $h,$ 
\[
S(x)\leq S(x-h)+S(h)\leq c(x-h)+S(h). 
\]%
Taking limits as $h\downarrow 0$ through such $h,$ yields%
\[
S(x)\leq cx, 
\]%
since $S(0+)=0$, by Theorem $0^{+}.$ This holds for any $x\in (0,1)$ and by
assumption for $x=0.$

For the last part, a similar appeal to subadditivity and Theorem $0^{+}$
yields $S(x+)\leq S(x),$ for any $x>0.$ So, for any $x\geq 0,$ if $x+h\in D,$
then%
\[
S(x)\geq \lim \sup_{h>0}S(x+h)\geq \lim \sup_{h>0,x+h\in D}S(x+h)=\lim
\sup_{h>0,x+h\in D}c(x+h)=cx. 
\]%
Combining, $S(x)=cx.$ $\square $

\bigskip

A similar but simpler argument yields:

\bigskip

\noindent \textbf{Proposition 3}$^{\prime }$\textbf{. }\textit{For} $S:%
\mathbb{R}\rightarrow \mathbb{R\cup \{-\infty },\mathbb{+\infty \}}$\textit{%
\ subadditive with }$S(0)=0$\textit{, write }$\Sigma _{+}:=\{u\in \mathbb{R}%
:S(u)<\infty \}.$ \textit{If }$\Sigma _{+}\cap \lbrack 0,\infty )$ \textit{%
contains an open set accumulating on the right at }$0,$\textit{\ and }$S$ 
\textit{is finite on a dense subset of }$\mathbb{R}_{+}$\textit{, then }$S|%
\mathbb{R}_{+}$\textit{\ is finite. }

\textit{If, further, }$S$\textit{\ is finite on a subset unbounded below},%
\textit{\ then also }$S|\mathbb{R}_{-}$ \textit{is finite.}

\bigskip

\noindent \textbf{Proof.} For the first assertion, argue just as before,
replacing `$\leq cx$' above by `$<\infty $' to show that $\Sigma _{+}\cap 
\mathbb{R}_{+}=\mathbb{R}_{+}$. Then $S(x)<\infty $ for all $x>0.$ For any $%
x>0$ choose $u_{0}>x$ with $S(u_{0})$ finite; then, as $u_{0}-x>0,$%
\[
S(x)\geq S(u_{0})-S(u_{0}-x)>-\infty . 
\]

For the last part, if $S$ is finite at $u_{0}<0,$ then for $u\in (u_{0},0)$%
\[
S(u)\leq S(u_{0})+S(x-u_{0})<\infty :\qquad \mathbb{R}_{-}\subseteq \Sigma
_{+}. 
\]
Finally, $S|\mathbb{R}_{-}$ is finite, since $-\infty <-S(u)\leq S(-u)$, for 
$u>0.$ $\square $

\bigskip

\noindent \textbf{Proposition 5.} \textit{For} $S:\mathbb{R}\rightarrow 
\mathbb{R}$\textit{\ subadditive:\newline
\noindent }(i) \textit{if }$S$\textit{\ is bounded above on some interval,
say by }$K$\textit{\ on }$B_{\delta }(a)$\textit{, for instance if} $(HS(S))$%
\textit{\ holds, then for any }$b\in \mathbb{R}$%
\[
S(b+a)-K\leq S(x)\leq S(b-a)+K\qquad (x\in B_{\delta }(b)), 
\]%
\textit{\ in particular it is locally bounded;\newline
}\noindent (ii) \textit{if }$S$\textit{\ is locally bounded, then }$\lim
\inf_{t\rightarrow 0}S(t)\geq 0,$ \textit{so }$S(0+)=0$\textit{\ if }$%
(HS(S)) $\textit{\ holds.}

\bigskip

\noindent \textbf{Proof.} (i) If $S$ is bounded above by $K$ on $B_{\delta
}(a)=a+(-\delta ,+\delta ),$ then for any $b\in \mathbb{R}$, $S\ $is bounded
on $b+(-\delta ,+\delta ).$ Indeed, for any $x\in b+(-\delta ,+\delta ),$
since both $a+(x-b)$ and $a-(x-b)$ are in $B_{\delta }(a),$%
\begin{equation}
\left. 
\begin{array}{c}
S(x)\leq S(b-a)+S(a+x-b)\leq S(b-a)+K, \\ 
S(x)\geq S(b+a)-S(a+b-x)\geq S(b+a)-K.%
\end{array}%
\right\}  \tag{$\star $}
\end{equation}%
If $S$ satisfies $(HS)$, then $S(x)<1$ for $x\in (0,\delta )$ for some $%
\delta >0.$

\noindent (ii) Following [HilP, 7.4.3], select a sequence $\{z_{n}\}$ with $%
z_{n}\rightarrow 0$ and $S(z_{n})\rightarrow \lambda
_{-}:=\inf_{u\rightarrow 0}S(u).$ By local boundedness, $\lambda _{-}$ is
finite, so for any $\varepsilon >0$ and $n$ large enough $\lambda
_{-}-\varepsilon \leq S(2z_{n})\leq 2S(z_{n})<2(\lambda _{-}+\varepsilon ).$
So $\lambda _{-}\leq 2\lambda _{-},$ yielding $\lambda _{-}\geq 0.$ $\square 
$

\bigskip

Having motivated right-sided continuity as in Theorem 0, we now prove it.

\bigskip

\noindent \textbf{Proof of Theorem 0. }The condition is evidently necessary.
As for sufficiency, suppose given $z_{n}\uparrow 0$ as in the hypothesis, $%
x_{n}\uparrow 0$ with $S(x_{n})\rightarrow \lambda ,$ and any $\varepsilon
>0;$ we will show that $\lambda =0.$

Choose $\delta >0$ with $S(t)\leq \varepsilon $ for $0\leq t\leq \delta .$
W.l.o.g. we assume that $z_{1}>-\delta .$ Now choose $m(n)$ for $n\in 
\mathbb{N}$ with $z_{n}\leq x_{m(n)}.$ Then $0\leq x_{m(n)}-z_{n}\leq \delta 
$, as $x_{m(n)},z_{n}\in (-\delta ,0),$ and so%
\[
S(x_{m(n)})\leq S(x_{m(n)}-z_{n})+S(z_{n})\leq \varepsilon +S(z_{n}). 
\]%
Passing to the limit gives%
\[
\lambda \leq \varepsilon +0=\varepsilon . 
\]%
Taking limits as $\varepsilon \downarrow 0$ gives $\lambda \leq 0.$ But, as $%
-x_{m(n)}\in (0,\delta ),$%
\[
0=S(0)\leq S(x_{m(n)})+S(-x_{m(n)})\leq S(x_{m(n)})+\varepsilon , 
\]%
so taking limits gives%
\[
0\leq \lambda +\varepsilon , 
\]%
so $\lambda \geq 0$, as above. Combining, $\lambda =0,$ so $S$ is continuous
at $0.$ Finally, if $S$ is finite-valued and continuous at $0$, it is so at
any $x,$ since%
\[
-S(-h)\leq S(x+h)-S(x)\leq S(h).\text{\qquad }\square 
\]

\bigskip

\noindent \textbf{Proof of Theorem 0}$^{\prime }$\textbf{.} Since $S|\Sigma $
is continuous at $0$ it is bounded above on $\Sigma _{\delta }:=\Sigma \cap
(-\delta ,\delta )$ for some $\delta >0;$ but $\Sigma _{\delta }+\Sigma
_{\delta }$ contains an interval, so $S$ is bounded on an interval, and so
locally bounded by Prop. 5(i). If $S$ is not continuous at $0$, then by
Prop. 5(ii) $\lambda _{+}:=\lim \sup_{t\rightarrow 0}S(t)>\lim
\inf_{t\rightarrow 0}S(t)\geq 0.$ Choose a null sequence $\{z_{n}\}$ with $%
S(z_{n})\rightarrow \lambda _{+}>0.$ Let $\varepsilon :=\lambda _{+}/4.$
W.l.o.g. $S(z_{n})>\lambda _{+}-\varepsilon $ for all $n.$ By continuity on $%
\Sigma $ at $0$ there is $\delta >0$ with $|S(t)|<\varepsilon $ for $t\in
\Sigma _{\delta }.$ As before and using symmetry, $\Sigma _{\delta }+\Sigma
_{\delta }=\Sigma _{\delta }-\Sigma _{\delta }$ contains an interval $I$
around $0.$ For any $n$ with $z_{n}\in I,$ there are $u_{n},v_{n}\in \Sigma
_{\delta }$ with $z_{n}=u_{n}+v_{n};$ then%
\[
S(z_{n})\leq S(u_{n})+S(v_{n})\leq 2\varepsilon <\lambda _{+}, 
\]%
and so%
\[
3\lambda _{+}/4=\lambda _{+}-\varepsilon <S(z_{n})\leq S(u_{n})+S(v_{n})\leq
2\varepsilon <\lambda _{+}/2, 
\]%
a contradiction. So $S$ is continuous at $0$ and so continuous everywhere,
as in Theorem 0. The last part follows since $\Sigma \cap (0,\delta ),$
being Baire/measurable and non-negligible, has the SW property for each $%
\delta >0$. $\square $

\bigskip

\noindent \textbf{Proof of Theorem 0}$^{+}$\textbf{.} We may take $c=0$,
since $S^{\prime }(t):=S(t)-ct$ is linearly bounded above by $0$ on $\Sigma
, $ and $S^{\prime }$ is subadditive. (Also the thinned Heiberg-Seneta
condition holds for $S^{\prime }.)$ Thus $S(t)\leq 0$ for $t\in \Sigma .$

Fix an arbitrary $x>0.$ We show that $S(x)\leq 0.$ As $\Sigma $ accumulates
at $0,$ there is a point $\sigma _{x}\in \Sigma \cap (0,x/2).$ Then $\Sigma
^{\prime }:=\Sigma \cap (\sigma _{x},\frac{1}{2}(\sigma _{x}+x/2))$ has the
SW property locally, and so $\Sigma ^{\prime }+\Sigma ^{\prime }$ contains a
proper interval $[a,b]$ in $(2\sigma _{x},\sigma _{x}+x/2).$

With $a,b$ fixed, choose $\sigma \in \Sigma \cap (0,b-a)\cap (0,a).$ By
density of $\Sigma ,$ we may suppose that $a,x\notin \mathbb{N}\sigma .$
Then there is $m\in \mathbb{N}$ with $a<m\sigma <m\sigma +\sigma <b.$ Now,
as $(m+1)\sigma <b<x,$ we may choose $n>m+1$ in $\mathbb{N}$ with $n\sigma
<x<n\sigma +\sigma .$ Then, as $0<x-n\sigma <\sigma ,$ adding $m\sigma $
gives%
\[
a<m\sigma <x+(m-n)\sigma <(m+1)\sigma <b. 
\]%
Now pick $u,v\in \Sigma ^{\prime }$ with $u+v\in (a,b)$ such that%
\[
u+v=x-(n-m)\sigma :\qquad x=u+v+(n-m)\sigma . 
\]%
By subadditivity, as $n-m\in \mathbb{N}$, and as $u,v,\sigma \in \Sigma ,$%
\[
S(x)\leq (n-m)S(\sigma )+S(u)+S(v)\leq 0. 
\]%
Thus $S(x)\leq cx$ for all $x>0.$ In particular%
\[
\lim \sup_{x\downarrow 0}\text{ }S(x)\leq 0. 
\]

Being linearly bounded above on $\Sigma $, $S$ is also relatively locally
bounded above on $\Sigma $, hence also on $\Sigma +\Sigma ,$ and so on an
interval; so by Prop. 5(ii) $\lim \inf_{x\downarrow 0}S(x)\geq 0.$

The final assertion follows from Theorem 0. $\square $

\bigskip

We may now turn to linearity from additivity rather than subadditivity, for
which see later. A key step follows. By appealing to Kronecker's Theorem
([HarW, XXIII, Th. 438]), the proof rolls together the two ingredients of
density and of routine use of continuity (as in Prop. 2$^{\prime }$ above);
we thank the Referee for this elegant approach. See Proposition $6^{\prime }$
in \S 5 for an alternative approach.

\bigskip

\noindent \textbf{Proposition 6. }\textit{For }$\mathbb{A}$ \textit{a dense
subgroup of }$\mathbb{R}$, \textit{if }$G:\mathbb{A\rightarrow R}$ \textit{%
is additive and bounded above on }$(0,\varepsilon )\cap \mathbb{A}$\textit{,
for some }$\varepsilon >0,$\textit{\ then }$G$\textit{\ is linear: }$G(a)=ca$%
\textit{\ for some }$c\in \mathbb{R}$ \textit{and all }$a\in \mathbb{A}$.

\bigskip

\noindent \textbf{Proof. }Being subadditive, $G$ may be assumed bounded on $%
I=(0,\varepsilon )$ for some $\varepsilon >0,$ by Prop. 5(i).

Fix any non-zero $u_{0}\in \mathbb{A}$ and put $c=G(u_{0})/u_{0}.$ We prove
that $G(a)=ca$ for all $a\in \mathbb{A}$ by showing that $H(a):=G(a)-ca$ ($%
a\in \mathbb{A}$) is identically zero.

Now $H$ is bounded on $I\cap \mathbb{A}$ , by $M$ say. By additivity, $%
H(a)=0 $ for $a\in u_{0}\mathbb{Z}$, as $H(u_{0})=0.$ Suppose that $H(u)\neq
0$ for some $u\in \mathbb{A}$; then for any $p\in \mathbb{N}$, $pu\notin
u_{0}\mathbb{Z}$ (as otherwise $pH(u)=H(pu)=0$). Fix $p\in \mathbb{N}$
arbitrarily. As $pu$ and $u_{0}$ are incommensurable, the subgroup they
generate, $u_{0}\mathbb{Z+}pu\mathbb{Z}$, is dense, by Kronecker's theorem.
So $mpu+nu_{0}\in I\cap (0,|u_{0}|)\cap \mathbb{A},$ for some $m,n\in 
\mathbb{Z}$. As $m\neq 0$ (since $nu_{0}\in (0,|u_{0}|)$ is not possible), 
\[
M\geq |H(mpu+nu_{0})|=|m|p\text{ }|H(u)|\geq p|H(u)|. 
\]%
So $M/p\geq |H(u)|.$ But $p\in \mathbb{N}$ was arbitrary, so $H(u)=0,$ after
all. $\square $

\bigskip

\noindent \textbf{Proposition 7 (Unique extension).} \textit{For }$\Sigma
\subseteq \mathbb{R}$ \textit{dense and closed under integer scaling (e.g. a
subgroup), let }$G:$\textit{\ }$\Sigma \rightarrow \mathbb{R}$ \textit{be
linear: }$G(\sigma )=c\sigma $\textit{\ (}$\sigma \in \Sigma $\textit{). If} 
$S:\mathbb{R}\rightarrow \mathbb{R}$ \textit{with }$S(0+)=S(0)$ \textit{is
any subadditive extension of }$G$\textit{, then }$S$\textit{\ is also linear
on }$\mathbb{R}$\textit{\ and }$S(u)=S(1)u=cu$\textit{\ for all }$u\in 
\mathbb{R}$\textit{.}

\bigskip

\noindent \textbf{Proof. }Here $S(0)=0,$ since $S(-\sigma )=G(-\sigma
)=-G(\sigma )=-S(-\sigma )$ for any $\sigma \in \Sigma .$ Take $z_{n}\in
\Sigma \cap (-1,0)$ converging to $0;$ then $S(z_{n})=G(z_{n})=cz_{n}%
\rightarrow 0$ for some $c.$ By Theorem 0, $S$ is continuous. Now $S(\sigma
)=G(\sigma )=c\sigma $ on $\Sigma ;$ so, as $\Sigma $ is dense in $\mathbb{R}
$, by continuity $S(t)=ct=tS(1)$ on $\mathbb{R}$. $\square $

\bigskip

Linearity from subadditivity is now a corollary of Prop. 5, 6 and 7:

\bigskip

\noindent \textbf{Proof of Theorem 1. }Here $S(0)=0$ as $S|\mathbb{A}$ is
additive on $\mathbb{A}$. As ($HS(S)$) holds, for each $\varepsilon >0$
there is $\delta >0$ with $S(t)\leq \varepsilon $ for $0<t<\delta $ and so $%
S(0+)=S(0).$ In particular $G=S|\mathbb{A}$ is additive on $\mathbb{A}$ and
bounded above on $(0,\delta ).$ By Prop. 6$,$ $G(a)=ca$ for some $c\in 
\mathbb{R}$ and all $a\in \mathbb{A}$. By Prop. 7, $S(t)=ct$ for all $t\in 
\mathbb{R}$. $\square $

\bigskip

\noindent \textbf{Proof of Theorem 1}$^{\prime }.$ \textbf{(a)} This follows
from Theorem 1 by Theorem 0$^{\prime }.$ $\square _{a}$

\noindent \textbf{(b)} This follows from Theorem 1 by Theorem 0$^{+}.$ $%
\square _{b}$

\noindent \textbf{(c) }Here $\mathbb{A\supseteq }\Sigma +\Sigma ,$ so
contains an interval, and, being a dense additive subgroup, $\mathbb{A}=%
\mathbb{R}$. So $S\ $is additive and is bounded on $\Sigma +\Sigma $, and so
on some interval; so $S$ is linear, by Prop. 6. $\square _{c}$

\noindent \textbf{(d) }As before, $S$ is subadditive, and by assumption is,
in view of Prop. 1, finite on the dense (additive) subgroup $\mathbb{A}=%
\mathbb{A}_{S}$ of $\mathbb{R}.$ As\textbf{\ }$\mathbb{A}$ is shift-compact, 
$\mathbb{A\supseteq A}+\mathbb{A}$ contains an interval, so by density again 
$\mathbb{A}=\mathbb{R}$, i.e. $S$ is finite everywhere. So $S=S|\mathbb{A}$
is \textit{additive} on the subgroup $\mathbb{A}=\mathbb{R}.$ In particular $%
[S(u+x)-S(x)]=S(u)$ and so $G(u)=S(u)$ for $u\in \mathbb{A}.$ By Prop. 5, $S$
is locally bounded, and so by Prop. 6 $S$ is linear on $\mathbb{A}=\mathbb{R}
$. $\square _{d}$ $\square $

\bigskip

\noindent \textbf{Theorem 2 (On linearity). }\textit{Let} $S:\mathbb{R}%
\rightarrow \mathbb{R}$ \textit{be subadditive with }$S(0+)=0$\textit{.} 
\textit{If, for some dense additive subgroup }$\mathbb{A}$ \textit{of }$%
\mathbb{R}$\textit{, the restriction }$S|\mathbb{A}$\textit{\ is additive,} 
\textit{then }$S$\textit{\ is linear on }$\mathbb{R}$\textit{: }$S(u)=cu$%
\textit{\ for some} $c\in \mathbb{R}$\textit{\ and all }$u\in \mathbb{R}$%
\textit{.}

\bigskip

\noindent \textbf{Proof.} $S$ extends $G:=S|\mathbb{A}$, additive and
continuous by Prop. 7. $\square $

\bigskip

\noindent \textbf{Cautionary Example. }Recall from \S 1 that $S:=\mathbf{1}_{%
\mathbb{R}\backslash \mathbb{Q}}$ is subadditive, and for $\mathbb{A}=%
\mathbb{Q}$, $S|\mathbb{A}$ $=0$ is linear; but $S$ is not linear. We return
to the relation of this example to Theorem 2 later in \S 7.

\bigskip

\noindent \textbf{Remark. }In Theorem 0$^{\prime }$ above, it is not enough
to assume only that the subadditive function $S\ $is locally bounded above
and $\mathbb{Q}$-homogeneous on a set $\Sigma $ that is dense (i.e. $%
S(q\sigma )=qS(\sigma )$ \textit{for all }$\sigma \in \Sigma $ \textit{and
rational }$q$); indeed, the indicator function $\mathbf{1}_{\mathbb{R}%
\backslash \mathbb{Q}}$ just considered is subadditive and also $\mathbb{Q}$%
-homogeneous on $\Sigma =\mathbb{Q}$, but not continuous. Such a weaker
assumption yields only that $\lim \inf_{t\rightarrow 0}S(t)=0=S(0).$\textit{%
\ }(Proof:\textit{\ }As $S$ is locally bounded, choose $K$ and $\delta >0$
with $S$ bounded by $K$ on $[-\delta ,\delta ].$ Fix $\sigma \in \Sigma \cap
(0,\delta );$ then $|S(\sigma /n)|=|S(\sigma )|/n\leq K/n$ for all $n\in 
\mathbb{N}$.)

\section{ Thinnings of Heiberg-Seneta conditions and subadditive functions}

We turn now from functional equations (whose prototype for us is the Cauchy
functional equation (Add) from \S 1) to functional inequalities (the
prototype of which for us is the corresponding inequality (Sub) from \S 1).
The classical sources here are [HilP, Ch. 3] (for the measurable case only,
but we need the category version also, for which see [BinO1]) and [Kuc, Ch.
16]; cf. [MatS]). Kuczma makes the contrast between the surprisingly great
affinity between Cauchy's equation and Jensen's inequality, and the
differences between (Add) and (Sub). Here matters are reversed: what is
surprising in our context is the \textit{similarity }between (Add) and (Sub).

The motivation in this section is the quantifier weakening of our title, in
the context of regular variation (\S 7 below). The prototypical results here
are \textit{characterization theorems} (BGT, Th. 1.4.3, 3.2.4). The
prototypical conditions for these are the \textit{Heiberg-Seneta conditions}%
, of `limsup limsup' type. There are links with \textit{Tauberian conditions}%
, of `lim limsup sup' type; see e.g. BGT Ch. 4. Here we start with the
classical condition ($HS(S)$) 
\[
\lim \sup\nolimits_{u\downarrow 0}S(u)\leq M 
\]%
for $M\in \mathbb{R}$, re-phrased as follows:%
\[
\limsup_{n\rightarrow \infty }\sup \{S(x):x\in (0,1/n)\}\leq M, 
\]%
which we generalize so as to `thin out' the intervals $(0,1/n)$ in various
senses (including with the help of category and measure), appropriately
expanding the notation $(\lim \sup )_{M}$. We focus here on sets that have
the local Steinhaus-Weil property of \S 1, and begin with a list of those
relevant here. For an alternative mode of thinning (via spanning) see \S 6.

\bigskip

\noindent \textbf{Examples} \textbf{of families of locally Steinhaus-Weil
sets}.

The sets listed below are typically, though not always, members of a
topology on an underlying set.\footnote{%
Below we refer to ideal topologies in the sense of [LukMZ].}

\noindent (o) $\Sigma $ a usual (Euclidean) open set in $\mathbb{R}$ (and in 
$\mathbb{R}^{n})$ -- this is the `trivial' example;

\noindent (i) $\Sigma $ density-open subset of $\mathbb{R}$ (similarly in $%
\mathbb{R}^{n})$ (by Steinhaus's Theorem -- see e.g. BGT Th. 1.1.1,
[BinO18], [Oxt, Ch. 8]);

\noindent (ii) $\Sigma $ Baire, locally non-meagre at all points $x\in
\Sigma $ (by the Piccard-Pettis Theorem -- as in BGT Th. 1.1.2, [BinO18],
[Oxt, Ch. 8] -- such sets can be `thinned out', i.e. extracted as subsets of
a second-category set, using separability or by reference to the Banach
Category Theorem [Oxt, Ch.16]);

\noindent (iii) $\Sigma $ the Cantor `middle-thirds excluded' subset of $%
[0,1]$ (since $\Sigma +\Sigma =[0,2]);$

\noindent (iv) $\Sigma $ universally measurable and open in the \textit{ideal%
} topology ([LukMZ], [BinO17]) generated by omitting Haar null sets (by the
Christensen-Solecki Interior-points Theorem of [Sol]);

\noindent (v) $\Sigma $ a Borel subset of a Polish abelian group and and
open in the ideal topology generated by omitting \textit{Haar meagre} sets
in the sense of Darji [Darj] (by Jab\l o\'{n}ska's generalization of the
Piccard Theorem, [Jab1, Th.2], cf. [Jab3], and since the Haar-meagre sets
form a $\sigma $-ideal [Darj, Th. 2.9]); for details see [BinO18].

If $\Sigma $ is \textit{Baire} (has the Baire property) and is locally
non-meagre, then it is co-meagre (since its quasi interior is everywhere
dense).

\noindent \textbf{Caveat. }1. Care is needed in identifying locally SW sets:
Mato\u{u}skov\'{a} and Zelen\'{y} [MatZ] show that in any non-locally
compact abelian Polish group there are closed non-Haar null sets $A,B$ such
that $A+B$ has empty interior. Recently, Jab\l o\'{n}ska [Jab4] has shown
that likewise in any non-locally compact abelian Polish group there are
closed non-Haar meager sets $A,B$ such that $A+B$ has empty interior.

\noindent 2. For an example on $\mathbb{R}$ of a compact subset $S$ such
that $S-S$ contains an interval, but $S+S$ has measure zero and so does not,
see [CrnGH].

\noindent 3. Below we are concerned with subsets $\Sigma \subseteq \mathbb{R}
$ where such `anomalies' are assumed not to occur.

\bigskip

\noindent \textbf{Definition. }Say that $(SW)$\textrm{-}$\lim
\sup_{u\downarrow 0}S(u)\leq M\in \mathbb{R}$ if there is $\Sigma \subseteq
(0,1)$ \textit{accumulating} at $0$ with the local Steinhaus-Weil property
such that, for $\Sigma _{n}:=\Sigma \cap (0,1/n),$%
\begin{equation}
\limsup_{n\rightarrow \infty }\sup S(\Sigma _{n}):=\limsup_{n\rightarrow
\infty }\sup \{S(x):x\in \Sigma _{n}\}\leq M.  \tag{$SW$-$\lim \sup_{M}(S)$}
\end{equation}

Evidently, $(SW)$\textrm{-}$\lim \sup_{u\downarrow 0}S(u)\leq M\in \mathbb{R}
$ holds if $\lim \sup_{u\downarrow 0}S(u)\leq M\in \mathbb{R}$ holds (refer
to $\Sigma =(0,1)).$

For later use, say that $(SSW)$\textrm{-}$\lim \sup_{u\downarrow 0}S(u)\leq
M\in \mathbb{R}$ if there is a \textit{symmetric} set $\Sigma \subseteq
(-1,1),$ i.e. $\Sigma =-\Sigma ,$ with the local Steinhaus-Weil property and 
$\Sigma \cap (0,1)$ accumulating at $0$, such that $SW$\textrm{-}$\lim
\sup_{M}(S)$ above holds for $\Sigma _{n}:=\Sigma \cap (-1/n,1/n).$

It is thematic for us that, inasmuch as they affect subadditive functions,
the quantifier weakening that thinness offers implies a level of
informativeness equal to that of the trivial example (o). We thank the
Referee for the following `bridging' result, clarifying the relationship
with $\lim \sup_{u\downarrow 0}S(u)\leq M\in \mathbb{R}$.

\bigskip

\noindent \textbf{Proposition 8}. \textit{For} $S:\mathbb{R}\rightarrow 
\mathbb{R}$ \textit{subadditive, }$SW$-\textrm{limsup}$_{M}(S)$ \textit{%
holds for some} $M$ \textit{iff: for some }$K,$%
\begin{equation}
S\text{ is bounded above by }K\text{ on some open }U\subseteq (0,\infty )%
\text{ with }0\in \bar{U}.  \tag{$\dag $}
\end{equation}

\noindent \textbf{Proof. }Suppose that $SW$-$\lim \sup_{M}(S)$ holds for
some $M.$ Then, with the notation above, for some infinite set $\mathbb{M}$ 
\[
S(u)<M+1\qquad (u\in \Sigma _{m},m\in \mathbb{M}). 
\]%
So for $x=u+v$ with $u,v\in \Sigma _{m}$ 
\[
S(x)<2(M+1). 
\]%
By the Steinhaus-Weil property $\Sigma _{m}+\Sigma _{m}$ contains an
interval $I_{m}$ in $(0,2/m);$ then $U:=\bigcup\nolimits_{m\in \mathbb{M}%
}I_{m}$ is open, $0\in \bar{U}$ and $S$ is bounded above by $K:=2(M+1)$ on $%
P.$

The converse is clear: given $U$ and $K$ as in $(\dag )$ above, take $\Sigma
:=U$; then $SW$-$\lim \sup_{M}(S)$ holds for $M=K.$ $\square $

\bigskip

The occurrence above of the infinite set $\mathbb{M}$ justifies a
combinatorial departure `beyond Lebesgue and Baire'. A wider combinatorial
characterization involving the embedding of a convergent subsequence (rather
than only of a null subsequence) may be obtained by reference to the level
sets of a function $S,$ defined by%
\[
H^{r},\text{ or }H^{r}(S):=\{x:|S(x)|<r\}. 
\]

\noindent \textbf{Proposition 9. }\textit{If the subadditive function }$S:%
\mathbb{R\rightarrow R}$ \textit{is Baire/measurable, then for every
convergent sequence }$\{u_{n}\}$\textit{\ with }$u_{n}\rightarrow u,$\textit{%
\ there exist }$k\in \mathbb{N}$, $t\in H^{k}$ \textit{and} $\mathbb{M}$ 
\textit{infinite with}%
\[
\{t+u_{m}:m\in \mathbb{M}\}\subseteq H^{k}. 
\]%
\textit{In particular, }$S$\textit{\ is locally bounded above, and so
locally bounded.}

\bigskip

\noindent \textbf{Proof.} We argue as in [BinO3]: since $\mathbb{R}%
=\bigcup\nolimits_{k\in \mathbb{N}}H^{k},$ there is $K\in \mathbb{N}$ with $%
H^{K}$ non-negligible Baire/measurable and $S|H^{K}$ is bounded by $K.$
Given a sequence $\{u_{n}\}$ with $u_{n}\rightarrow u,$ put $%
z_{n}:=u_{n}-u\rightarrow 0;$ by the \textit{bilateral} version of KBD
[BinO5, \S 3], choose $t\in H^{K}$ and $\mathbb{M}$ infinite with $%
\{t+z_{m},t-z_{m}:m\in \mathbb{M}\}\subseteq H^{K}.$ Then, with $a=t,b=t+u$
and $x=t+u+z_{m}$ for $m\in \mathbb{M}$, we may apply $(\star )$ as in the
proof of Prop. 5 but with $H^{K}$ for $B_{\delta }(a)$, since $%
a+b-x=t-z_{m}\in H^{K}$ and $x+a-b=t+z_{m}\in H^{K}.$ This yields, since $%
b+a=u+2t$ and $b-a=u,$%
\[
S(2t+u)-K\leq S(t+u_{m})\leq S(u)+K. 
\]%
So for some $k$ and infinite $\mathbb{M}$ one has $\{t+u_{m}:m\in \mathbb{M}%
\}\subseteq H^{k}.$

Suppose that $S$ is not locally bounded above; then for some $u$ and
sequence $\{u_{n}\}$ with $u_{n}\rightarrow u,$ the sequence $\{S(u_{n})\}$
is unbounded above. But for some $k\in \mathbb{N}$, $t\in H^{k}$ and $%
\mathbb{M}$ infinite\textit{\ }as above, since $t+u_{m}\in H^{k}$ for $m\in 
\mathbb{M}$,%
\[
S(u_{n})\leq S(t+u_{m})+S(-t)\leq k+S(-t), 
\]%
a contradiction. Now apply Prop. 5(i). $\square $

\bigskip

\noindent \textbf{Remark.} For $\Sigma $ shift-compact and $S:\mathbb{R}%
\rightarrow \mathbb{R}$ subadditive, if $S\ $is bounded above on $\Sigma ,$
then $S$ is bounded above on $\mathbb{R}$ [BinO1, Th. 2(ii)]. Proposition 3$%
^{\prime \prime }$ below extends this to the range $\mathbb{R\cup \{-\infty }%
,\mathbb{+\infty \}}$.

\bigskip

\noindent \textbf{Proposition 3}$^{\prime \prime }$\textbf{\ (On finiteness).%
} \textit{Suppose that} $S:\mathbb{R}\rightarrow \mathbb{R\cup \{-\infty },%
\mathbb{+\infty \}}$ \textit{is subadditive and satisfies:\newline
}\noindent (i)\textit{\ }$S$ \textit{is finite on a subset }$\Sigma $ 
\textit{unbounded below} (\textit{e.g. a dense subset of }$\mathbb{R}$), 
\textit{and\newline
}\noindent (ii)\textit{\ for some }$M\in \mathbb{R},$ $SW$\textrm{-}$\lim
\sup {}_{M}(S)$ \textit{holds}$.$\textit{\newline
Then }$S$\textit{\ is finite everywhere, and so locally bounded.}

\bigskip

\noindent \textbf{Proof. }By (ii) $S$ is bounded above, by $K$ say, on some
locally-SW set $T$ (e.g. a Baire/measurable non-negligible); as above, by
the Steinhaus property, $T+T\ $contains an interval on which $S$ is bounded
above by $2K.$ Apply Proposition 3 to deduce that $S$ is finite everywhere.

As $S$ is finite and subadditive and bounded above on an interval, by Prop.
4 (cf. [BinO1, Th. 2(ii)]) $S$ is locally bounded. $\square $

\bigskip

We close with an alternative approach to Theorem $0^{\prime }.$

\bigskip

\noindent \textbf{Proposition 4}$^{\prime }$\textbf{.} \textit{If }$S:%
\mathbb{R}\rightarrow \mathbb{R}$\textit{\ is subadditive} \textit{with} $%
S(0)=0$ \textit{and satisfies }$SSW$-$\lim \sup (S)_{0}$\textit{, then }$S$ 
\textit{is continuous.}

\bigskip

\noindent \textbf{Proof.} Otherwise, $S$ is not continuous at $0$ (as at the
end of Prop. 11), and so as above $\lambda _{+}:=\lim \sup_{t\rightarrow
0}S(t)>\lim \inf_{t\rightarrow 0}S(t)\geq 0$, the latter by Prop. 5(ii) (by
the Remark above, combined with Prop. 5(i), $S$ is locally bounded). Choose
a null sequence $\{z_{n}\}$ with $S(z_{n})\rightarrow \lambda _{+}>0,$ and $%
\Sigma $ as in $SSW$-$\lim \sup (S)_{0}.$ For any $\varepsilon >0,$ take $%
N=N(\varepsilon )$ such that $\sup S(\Sigma _{N})<\varepsilon /4.$ As $%
\Sigma _{N}+\Sigma _{N}$ contains an open interval, $r+I$ with $0\in I$ say,
then $r+z_{n}\in r+I$ for all large $n;$ write $r=s+t$ and $%
r+z_{n}=u_{n}+v_{n}$ for $s,t,u_{n},v_{n}\in \Sigma _{N}$. Then 
\[
S(z_{n})\leq S(u_{n})+S(v_{n})+S(-s)+S(-t)<\varepsilon /4+\varepsilon
/4+\varepsilon /4+\varepsilon /4\leq \varepsilon , 
\]%
as $-s,-t\in \Sigma _{N}$ (symmetry). Taking limits yields $\lambda _{+}\leq
\varepsilon ,$ for each $\varepsilon >0,$ so $\lambda _{+}\leq 0,$ a
contradiction. So $\lambda _{+}=0,$ and so $\lim_{t\rightarrow
0}S(t)=0=S(0), $ contradicting our initial assumption$.$ $\square $

\section{Sublinearity and Berz's Theorem}

Theorem 2 above is reminiscent of the following classical result. Recall
that for $\Sigma $ closed under positive rational scaling $S$ is \textit{%
sublinear} on $\Sigma $ in the sense of Berz [Berz] if $S$ is subadditive
and $S(nx)=nS(x)$ for $x\in \Sigma ,n=0,1,2,...$ (i.e. $S$ is positively $%
\mathbb{Q}$-homogeneous and $S(0)=0)$.

\bigskip

\noindent \textbf{Theorem B (Berz's Theorem, }[Berz]; cf. [Kuc, Th. 16.4.3])%
\textbf{.} \textit{If }$S:\mathbb{R}\rightarrow \mathbb{R}$ \textit{is
measurable and sublinear, then there are }$c_{\pm }\in \mathbb{R}$\textit{\
such that }$S(u)=c_{+}u$ \textit{for }$u\geq 0$ \textit{and }$S(v)=c_{-}v$ 
\textit{for} $v\leq 0.$

\bigskip

Theorem 3 below is its category analogue; cf. [BinO16] and \S 8.4. We will
need the following variant of Prop. 6, which in fact includes it (below take 
$\Sigma =\mathbb{A}$ a dense subgroup of $\mathbb{R}$, and $S=A$ an additive
function: continuity below is then equivalent to local boundedness at $0,$
by Prop. 2$^{\prime }$).

\bigskip

\noindent \textbf{Proposition 6}$^{\prime }$\textbf{\ (Relative linearity). }%
\textit{Let} $S:\Sigma \rightarrow \mathbb{R}$ \textit{be subadditive with }$%
\Sigma \cap \mathbb{R}_{+}$\textit{\ dense on }$\mathbb{R}_{+}$\textit{\ and
closed under positive-integer scaling. If }$S$\textit{\ is }$\mathbb{N}$%
\textit{-homogeneous and right-continuous,} \textit{then }$S$\textit{\ is
linear on }$\Sigma \cap \mathbb{R}_{+}$\textit{: }$S(\sigma )=c_{+}\sigma $%
\textit{\ for some} $c_{+}\in \mathbb{R}$\textit{\ and all }$\sigma \in
\Sigma \cap \mathbb{R}_{+}.$

\bigskip

\noindent \textbf{Proof. }The proof is adapted from BGT Th. 3.2.5 (see also
[BinG, Proof of Th. 5.7]). Fix any \textit{positive} $\sigma _{0}\in \Sigma $
and put $c:=S(\sigma _{0})/\sigma _{0}$. We show that $S(\sigma )=c\sigma $
for all $\sigma \in \Sigma .$ To this end, fix any $\sigma \in \Sigma .$ Now
define for $0<\delta \in \Sigma $%
\[
i=i(\delta ):=\min \{n\in \mathbb{N}:n\delta >\sigma \},\text{ }%
i_{0}=i_{0}(\delta ):=\min \{n\in \mathbb{N}:n\delta >\sigma _{0}\}, 
\]%
so that $i(\delta )\delta \downarrow \sigma $ and $i_{0}(\delta )\delta
\downarrow \sigma _{0}$ as $\delta \downarrow 0;$ here w.l.o.g. $G(\delta
)\neq 0$ (otherwise $G=0$ on $(0,\varepsilon )\cap \Sigma $, for some $%
\varepsilon >0$, implying below that $S(\sigma )=0$ and so $S\equiv 0$ on $%
\Sigma \cap \mathbb{R}_{+}$).

Since $\delta ,u_{0}\in \Sigma \cap \mathbb{R}_{+}$ and $\Sigma \cap \mathbb{%
R}_{+}$ is closed under integer scaling, we have by $\mathbb{N}$-homogeneity
of $S$ that $S(i\delta )=iS(\delta )$; likewise $S(i_{0}\delta
)=i_{0}S(\delta )$. Taking limits here and below through $\Sigma $ as $%
\delta \downarrow 0$ and using right-continuity of $S$ at $\sigma _{0}$ and $%
\sigma ,$ 
\[
S(\sigma _{0})=\lim\nolimits_{\delta \downarrow 0}S(i_{0}(\delta )\delta ),%
\text{ }S(\sigma )=\lim\nolimits_{\delta \downarrow 0}S(i(\delta )\delta ). 
\]%
Dividing these two, as $\sigma _{0}\neq 0$,%
\[
\sigma /\sigma _{0}=\lim_{\delta \downarrow 0}i(\delta )\delta /i_{0}(\delta
)\delta =\lim_{\delta \downarrow 0}i(\delta )/i_{0}(\delta )=\lim_{\delta
\downarrow 0}i(\delta )S(\delta )/i_{0}(\delta )S(\delta )=S(\sigma
)/S(\sigma _{0}). 
\]%
Cross-multiplying, $S(\sigma )=c\sigma $ for any $u\in \Sigma \cap \mathbb{R}%
_{+}.$ $\square $

\bigskip

\noindent \textbf{Theorem 3 (Baire-Berz Theorem).} \textit{If }$S:\mathbb{R}%
\rightarrow \mathbb{R}$ \textit{is Baire and sublinear, there are }$c_{\pm
}\in \mathbb{R}$\textit{\ such that }$S(u)=c_{+}u$ \textit{for }$u\geq 0$ 
\textit{and }$S(v)=c_{-}v$ \textit{for} $v\leq 0.$

\bigskip

\noindent \textbf{Proof.} By a theorem of Baire, we may choose a meagre set $%
M$ such that $S|(\mathbb{R}\backslash M)$ is continuous (see [Oxt, Th.
8.1]). Expand $M$ to a union of closed nowhere dense sets, if necessary.
Take $\Sigma :=\mathbb{R}\backslash \bigcup\nolimits_{q\in \mathbb{Q}}qM$,
which is closed under rational scaling and is dense on $\mathbb{R}$ (by
Baire's Theorem -- as each $\mathbb{R}\backslash qM$ is a dense $\mathcal{G}%
_{\delta }$). In particular $\Sigma =-\Sigma .$ By Prop. 6$^{\prime }$, $S$
is linear on $\Sigma \cap \mathbb{R}_{+}$. By Theorem 0$^{+}$, $S$ is
continuous. So as $S$ is linear on $\Sigma \cap \mathbb{R}_{+},$ a dense
subset of $\mathbb{R}_{+}$, it is linear on $\mathbb{R}_{+}$; likewise $%
S(-x) $ is linear on $\mathbb{R}_{+}$. $\square $

\section{Thinning via spanning}

Weakening of quantifiers amounts to \textit{thinning} of the relevant set.
Theorem 1 may be viewed as a two-pronged test of the subadditive function $S$
on thin sets: one prong a domain condition, density, the other a boundedness
condition at $0$. As $\mathbb{A}$ is a subgroup, the first condition reduces
to density of $\mathbb{A}$ at $0.$ On $\mathbb{R}$ it further reduces to a
two-point condition: for Theorem 1 to hold, the set $\mathbb{A}$ necessarily
has at least two incommensurable members (recall the example $\mathbf{1}_{%
\mathbb{R}\backslash \mathbb{Q}})$; by Kronecker's Theorem ([HarW] XXIII,
Th. 438) $\mathbb{A}$ is then dense. That is, $\mathbb{A}$ needs to be at
least a two-dimensional subspace of $\mathbb{R}$, regarded as a vector space
over $\mathbb{Q}$.

Say that $T\subseteq \mathbb{R}$ is a \textit{spanning} set if it spans $%
\mathbb{R}$ regarded as a vector space over $\mathbb{Q}$ (e.g. contains a
Hamel basis). We work below with \textit{analytic sets}; for background see
e.g. [Rog], [BinO6]. The approach below is motivated by a theorem due to F.
Burton Jones and its later strengthening by Z. Kominek -- see below. Recall
from \S 4 that $T$ is \textit{symmetric} if $-T:=\{-t:t\in T\}=T,$ and that $%
T\ $is \textit{shift-symmetric} if for some $\tau $ the set $T+\tau
:=\{t+\tau :t\in T\}$ is symmetric ($T+\tau =-T-\tau )$, and the latter is
equivalent to \textit{self-similarity}, in the sense of [BinO9]: that $T=a-T$
for some $a$ (as then $T-a/2=-(T-a/2)$).

The vectorial view brings further insights. Jones [Jon] proved that for
additive $A:\mathbb{R}\rightarrow \mathbb{R}$, if $A|T$ is continuous on an
analytic spanning set $T$, then it is continuous. Much later, Kominek [Kom]
showed that, for such a set $T,$ if $A|T$ is bounded above on $T$, then $T$
is continuous -- cf. [BinO6, Th. JK] and \S 8.6. Kominek's Theorem is
stronger: it implies Jones's, for which see again [BinO6, \S 3]. The
Jones-Kominek results show \textit{how to test on thin sets }$T$\textit{\
properties of interest on} $\mathbb{R}.$

In Theorem 1$_{d}^{\prime }$ above the oscillation condition is further
thinned by using only shift-compact subsets in neighbourhoods of the origin%
\footnote{%
Say that $T$ is $k$\textit{-thinned} shift-compact if the $k$-fold sum: $%
k\cdot T:=T+...T$ is shift-compact. These thinned subsets could in principle
do the work of shift-compactness, as in [BinO4, \S 4]; the standard Cantor
set $C$ is $2$-thinned shift-compact, as $C+C\ $contains an interval (see
Remarks, below): boundedness of a subadditive function on $C$ yields local
boundedness.}. This raises the question of establishing further quantifier
weakening by testing only on analytic spanning sets, as these need not be
shift-compact. This is indeed possible on both prongs, as follows.

As we have seen, Theorem 1 rests on the presence of enough of the following
three properties: finiteness, local boundedness and continuity. These are
aided by density, and easy to achieve for additive functions by demanding
them on a spanning set. Kominek's Theorem cannot be applied directly to $S,$
as $S$ is only subadditive, but a simple modification of its proof will work
below. It is convenient to use a result of Erd\H{o}s, for which we need the
following.

\bigskip

\noindent \textbf{Definition.} For $H$ a Hamel basis, and by extension for $%
H $ a spanning set, we write%
\[
Z(H):=\{k_{1}h_{1}+...+k_{n}h_{n}:k_{i}\in \mathbb{Z},h_{i}\in H\} 
\]%
for the \textit{Erd\H{o}s set} of $H$ (for which see e.g. [Kuc, \S 11.5]).

\bigskip

\noindent \textbf{Theorem E (}Erd\H{o}s\textbf{; }see e.g. [Kuc, Lemma
11.5.3])\textbf{.} \textit{For }$H$\textit{\ a Hamel basis and, more
generally, for }$H$\textit{\ a spanning set, the Erd\H{o}s set }$Z(H)$ 
\textit{is a dense subgroup.}

\bigskip

We need a modified version of the Analytic Dichotomy Lemma in [BinO6, \S 2],
which refers to the $k$\textit{-fold sum} $k\cdot T:=T+...+T$ ($k$ times).

\bigskip

\noindent \textbf{Proposition 10.} \textit{If }$T\subseteq \mathbb{R}$%
\textit{\ is a symmetric analytic spanning set, then for some }$k\in \mathbb{%
N}$ \textit{the }$k$-\textit{fold sum }$k\cdot T$\textit{\ contains the
interval }$(-1,1).$ \textit{If }$T$\textit{\ is a shift-symmetric analytic
spanning set, then for some }$k\in \mathbb{N}$ \textit{the }$k$-\textit{fold
sum }$k\cdot T$\textit{\ contains an interval.}

\bigskip

\noindent \textbf{Proof.} We indicate the necessary modification to the
proof in [BinO6, \S 2]. As there, with $T$ an analytic spanning set, for
some $m$-tuple of rationals $r_{i}=p_{i}/N$ with $p_{i}\in \mathbb{Z}$, for $%
1\leq i\leq m,$ the set $(r_{1}T+...+r_{m}T)$ is non-null. So too is $%
(p_{1}T+...+p_{m}T),$ and so $(p_{1}T+...+p_{m}T)-(p_{1}T+...+p_{m}T)$
contains an interval around $0$ (by the Steinhaus property, again). So, in
the $T$ symmetric case, as $T=-T$ we conclude that $T^{\prime
}:=|p_{1}|T+...+|p_{m}|T+|p_{1}|T+...+|p_{m}|T$ contains an interval around $%
0,$ say $(-1/n,1/n).$ Since $n\cdot T^{\prime }\supseteq (-1,1),$ one has $%
k\cdot T\supseteq (-1,1)$ for $k=2n(|p_{1}|+...+|p_{m}|).$

If, however, $\tau +T\ $rather than $T$ is symmetric, then $kT=k(\tau
+T)-k\tau $ and the shift-symmetric case follows from the symmetric case. $%
\square $

\bigskip

We may now give the analytic-spanning version of Theorem 1$.$

\bigskip

\noindent \textbf{Theorem 4. }\textit{For} $S:\mathbb{R}\rightarrow \mathbb{%
R\cup \{-\infty },\mathbb{+\infty \}}$ \textit{subadditive, and }$\mathbb{A}%
\subseteq \mathbb{R}$\textit{\ \textit{an additive subgroup}, suppose that}%
\newline
\noindent (i) $S|\mathbb{A}$ \textit{is finite and additive;}\newline
\noindent (ii)\textit{\ there exists a shift-symmetric (self-similar),
analytic, spanning set }$T\subseteq \mathbb{A}\ $\textit{such that }$S|T$%
\textit{\ is locally bounded above.}\newline
\noindent \textit{Then }$S$\textit{\ is linear: }$S(u)=cu,$\textit{\ for some%
} $c\in \mathbb{R}$\textit{\ and all }$u\in \mathbb{R}.$

\bigskip

\noindent \textbf{Proof. }By Prop. 10, fix $k\in \mathbb{N}$ such that $%
k\cdot T$ contains $(-1,1).$ By Theorem E above, the Erd\H{o}s set $Z(T),$
and so also $\mathbb{A},$ is a dense subgroup. Now $k\cdot T\subseteq 
\mathbb{A}$, as $\mathbb{A}$ is a subgroup. As $k\cdot T$ contains an
interval, $\mathbb{A}=\mathbb{R}$ (cf. Th. S in \S 8.3). As $S|T$ is locally
bounded above, by subadditivity $S$ is locally bounded above on $k\cdot T$
(by continuity of addition), and so on an interval. But $S$ is additive on $%
\mathbb{R}$ and bounded above on an interval, so locally bounded, by Prop.
5(i). So $S$ is linear by Prop. 7. $\square $

\bigskip

The condition that $S$ be linearly bounded on a symmetric analytic spanning
set, stronger than (ii), yields automatic linearity -- without any need for
(i) above. There are echos here of Theorem 0$^{+}.$

\bigskip

\noindent \textbf{Proposition 11.}\textit{\ If the subadditive function }$S:%
\mathbb{R\rightarrow R}\ $\textit{is linear on a symmetric analytic spanning
set, then }$S$\textit{\ is linear.}

\bigskip

\noindent \textbf{Proof.} Suppose that $S(t)=ct$ for $t\in T$ with $T$ a
symmetric analytic spanning set. By Prop. 10 there exists $k\in \mathbb{N}$
with $k\cdot T$ containing $\,I:=(-1,1).$ So $\mathbb{R=}\dbigcup%
\nolimits_{n=k}^{\infty }n\cdot T,$ since $n\cdot T$ contains the interval $%
(n/k)I$ for $n\in \mathbb{N}$ with $n>k.$ Now consider any $k\in \mathbb{N}$
with $k\cdot T$ containing a symmetric interval $J$ about $0.$ For $u\in $ $%
J\subseteq kT,$ choose $t_{i}\in T$ with $u=\dsum\nolimits_{i=1}^{k}t_{i};$
then%
\[
S(u)\leq \dsum\nolimits_{i=1}^{k}S(t_{i})=\dsum\nolimits_{i=1}^{k}ct_{i}=cu, 
\]%
by linearity of $S$ on $T.$ Likewise, as $-u\in J,$%
\[
S(-u)\leq c(-u)=-cu. 
\]%
Note that $S(0)=0$, by subadditivity, as $S(-t)=-ct=-S(t)$ for $t\in T$; so $%
-S(u)\leq S(-u),$ again by subadditivity. Combining,%
\[
-cu\leq -S(u)\leq S(-u)\leq -cu:\qquad S(u)=cu\qquad (u\in I). 
\]%
As $J$ is symmetric and can have arbitrary length, $S(u)=cu$ for all $u$. $%
\square $

\bigskip

\noindent \textbf{Remarks. }1.\textbf{\ }We raise the question of whether
symmetry can be omitted in Theorem 4.

\noindent 2. With $C$ the standard Cantor set, consider its translate $T:=C-%
\frac{1}{2},$ which is symmetric. It is a spanning set as $C+C=[0,2],$ $%
T+T=[-1,1],$ as in Proposition 10; so, assuming the Axiom of Choice \textrm{%
AC}, a Hamel basis $H$ may also be selected in $T,$ and $H$ is non-dense as $%
C$ and $T$ are. So a (nowhere dense) shift-symmetric Cantor set, such as $T,$
would suffice to test for local boundedness. Though density is not
explicitly mentioned, it remains present implicitly as $T,$ being
uncountable, contains incommensurables.

\noindent 3. In view of the Jones-Kominek theorems above, the question
arises as to whether a subadditive $S$ with $S|T$ right-continuous on an
analytic spanning set $T$ is right-continuous on a dense set. Note that $S:=%
\mathbf{1}_{\mathbb{R}\backslash \mathbb{Q}}$ is continuous on the
(analytic) spanning set $T:=\mathbb{R}\backslash \mathbb{Q}$, yet $S$ is not
continuous.

A continuous function $f:T\rightarrow \mathbb{R}$ is termed (\textit{%
extensibly continuous} or just) \textit{precompact} in [BinO6, \S 5] if $%
\{f(t_{n})\}$ is a Cauchy sequence whenever $\{t_{n}\}$ is a Cauchy sequence
in $T$ (cf. [Ful]). For $f:\mathbb{R}\rightarrow \mathbb{R}$ additive,
precompactness of $f|T$ on an analytic spanning set $T$ implies continuity,
but this feature does not extend to subadditive functions, as $\mathbf{1}_{%
\mathbb{R}\backslash \mathbb{Q}}|(\mathbb{R}\backslash \mathbb{Q)}$ is
precompact.

\bigskip

\noindent \textbf{Definition. }Altering the property $SW$-$\lim \sup (S)_{0}$
so that the limsup is taken with reference to symmetric sets $P_{n}\subseteq
(-1/n,1/n)$ of the scaled form $P_{n}:=T/n=\{t/n:t\in T\},$ with $T$ a fixed 
\textit{symmetric} analytic spanning set in $(-1,1)$, yields $JK(T)_{0}$ --
the $(JK)_{0}$ \textit{property for} $T.$

\bigskip

This yields an analogue of Theorem 4 for general subadditive functions $S$
rather than those that are linear on an additive subgroup $\mathbb{A}$
containing an analytic spanning set $T$ with $S|T$ bounded. As each $P_{n}$
is a spanning set, the semigroup argument underpinning Prop. 3$^{\prime }$
remains valid, and so Prop. 3$^{\prime \prime }$ holds with $(JK)_{0}$
replacing $(SW$-$\lim \sup )_{M}.$

\bigskip

\noindent \textbf{Theorem 5 (Automatic right-continuity, after Jones-Kominek)%
}. \textit{If} $S:\mathbb{R}\rightarrow \mathbb{R}$ \textit{is subadditive,
and }$JK(T)_{0}$ \textit{holds for some symmetric analytic spanning set} $T$%
\textit{\ with }$S|T$ \textit{locally bounded above --} \textit{then }$%
S(0+)=0.$

\bigskip

\noindent \textbf{Proof.} By the Compact Spanning Approximation Theorem of
[BinO6, \S 3], passing to a symmetric compact subset of $T$ if necessary, we
may assume that $T\ $is compact, and, scaling if necessary, that $T\subseteq
(-1,1)$. By Prop. 10, fix $k\in \mathbb{N}$ such that $k\cdot T$ contains $%
(-1,1).$ If $S$ is locally bounded above on $T$ by $M,$ then $S$ is locally
bounded above on $(-1,1)$ by $kN,$ as $k\cdot T$ contains $(-1,1)$; so $S$
is locally bounded above, and so locally bounded by Prop. 5.

With this in mind, consider an arbitrary sequence $\{v_{n}\}$ with $%
v_{n}\rightarrow 0$ and $S(v_{n})\rightarrow a.$ For each $n,$ as $k\cdot
(T/n)$ contains $(-1/n,1/n),$ we may assume, by passing to a subsequence of $%
\{v_{n}\}$ if necessary, that $|v_{n}|\leq 1/n,$ and so $%
v_{n}=u_{n}^{1}+...+u_{n}^{k}$ with $u_{n}^{i}\in T/n$ for each $%
i=1,2,...,k. $ Again by passing to a subsequence if necessary, by $JK_{0}(T)$
we may assume that $\lim_{n}S(u_{n}^{i})\leq 0$ for $i=1,...,k.$ By
subadditivity, 
\[
S(v_{n})\leq S(u_{n}^{1})+...+S(u_{n}^{k}), 
\]%
so $a=\lim_{n}S(v_{n})\leq 0.$ Suppose that $a<0.$ W. l. o. g., $%
w:=\sum_{i}v_{i}<\infty .$ Put $w_{n}:=\sum_{i\leq n}v_{i}\rightarrow w$. By
subadditivity,%
\[
S(w_{n})\leq \sum\nolimits_{i\leq n}S(v_{i})\rightarrow -\infty . 
\]%
This contradicts local boundedness at $w$. So $a=0;$ i.e. $S(0+)=0.$ $%
\square $

\bigskip

\noindent \textbf{Remark.} The very last step in the proof above is inspired
by the Goldie argument in BGT, p. 142.

\section{Quantifier weakening in regular variation}

The standard work on the Karamata theory of regular variation is BGT. The
present authors have returned to this area in a number of papers, together
and separately, largely addressed to matters left open there. First, we
address the \textit{foundational question}: what is the appropriate
generalization of the measurability and Baire-property settings of BGT?
Secondly, we address the \textit{contextual question}: what, beyond the
real-line setting of BGT (and other settings briefly addressed in BGT
Appendix 1 such as the complex plane, Euclidean space and topological
groups), is the natural context for the theory? In Theorem 6 here, we
complete our reduction of the number of hard proofs in the area to zero,
thus making good on a claim we have already made elsewhere (see \S 1). It is
striking that Th. 1.4.3 of BGT, in the context of Karamata Theory,%
\begin{equation}
\lim\nolimits_{x\rightarrow \infty }f(\lambda x)/f(x)=g(\lambda )\text{%
\qquad }(\forall \lambda >0),  \tag{$K$}
\end{equation}%
is no harder than Th. 3.2.5 of BGT, in the context of \textit{Bojani\'{c}%
-Karamata/de Haan} Theory ([BojK], [dHa]; cf. [BinO14])%
\begin{equation}
\lim\nolimits_{x\rightarrow \infty }\{f(\lambda x)-f(x)\}/g(x)=k(\lambda )%
\text{\qquad }(\forall \lambda >0).  \tag{$BKdH$}
\end{equation}

Here we weaken the quantifier $\forall $ above as much as possible (cf.
[BinG]). Results of this nature are false if the quantifier is weakened too
much, for reasons connected with Hamel pathology (BGT\ \S\ 1.1.4; cf. [Kuc,
Ch. 11]). As is usual in this area, of course, one encounters a dichotomy:
matters are either very nice or very nasty -- even the merest hint of good
behaviour being sufficient to guarantee the former; cf. [BinO7, 11]. Here,
the condition needed, the \textit{Heiberg-Seneta condition}, is a \textit{%
one-sided }one of `liminf liminf' type, as in the BGT\ results cited above,
[BinG], [Hei], [Sen], or equivalently (below) of `limsup limsup' type.

As usual for proofs, we work with Karamata theory written additively rather
than multiplicatively. As in Prop. 1 above, we write%
\[
G(u):=\lim\nolimits_{x\rightarrow \infty }F(u+x)-F(x) 
\]%
for the limit on the right where this exists; $\mathbb{A}_{F}$ the set on
which the limit exists; and $F^{\ast }$ for the limsup, again as in \S 1.

\bigskip

\noindent \textbf{Theorem 6 (Quantifier-Weakening Theorem, }cf. [BinO15, Th.
6]\textbf{).} \textit{With }$F^{\ast }$ \textit{and }$\mathbb{A}_{F}$\textit{%
\ as above, suppose that}\newline
\noindent (i)\textit{\ }$\mathbb{A}_{F}$ \textit{is dense in }$\mathbb{R}$%
\textit{;}\newline
\noindent (ii)\textit{\ }$F^{\ast }$\textit{\ satisfies the one-sided
Heiberg-Seneta boundedness condition }%
\begin{equation}
\lim \sup\nolimits_{u\downarrow 0}F^{\ast }(u)\leq 0  \tag{$HS$}
\end{equation}%
-- \textit{then }$\mathbb{A}_{F}=\mathbb{R}$\textit{\ and }$F^{\ast }$%
\textit{\ is linear: }$F^{\ast }(u)=\lim_{x\rightarrow \infty
}[F(u+x)-F(x)]=cu$\textit{\ for some} $c\in \mathbb{R}$\textit{, and all }$%
u\in \mathbb{R}.$

\bigskip

\noindent \textbf{Proof of Theorem 6. }By Prop. 1, 3 and 6, $F^{\ast }$ is a
finite, subadditive, right-continuous extension of $G$. So $G$ is continuous
on $\mathbb{A}_{F},$ and so linear by Prop. 6: $G(\sigma )=c\sigma $ for
some $c$ and all $\sigma \in \mathbb{A}_{F}.$ As $\mathbb{A}_{F}$ is dense,
by Prop. 7, $F^{\ast }(u)=cu$ for all $u.$ By Prop. 1, $\mathbb{A}_{F}=%
\mathbb{R}$ and $F^{\ast }(u)=G(u).$ $\square $

\bigskip

\noindent \textbf{Cautionary Example. }For $F:=\mathbf{1}_{\mathbb{R}%
\backslash \mathbb{Q}}$ and $q\in \mathbb{Q}$, one has $F(q+x)-F(x)=0$ for
all $x,$ and so $\mathbb{A}_{F}$ is dense. Also $F^{\ast }=\mathbf{1}_{%
\mathbb{R}\backslash \mathbb{Q}}.$ Indeed, fix $t\notin \mathbb{Q}$; then $%
F(t+q)-F(q)=1$ for $q\in \mathbb{Q}$ and $F(t+x)-F(x)\in \{-1,0\}$ for $%
x\notin \mathbb{Q}$, so that $F^{\ast }(t)=1$ and $\mathbb{A}_{F}=\mathbb{Q}$%
. As $F^{\ast }$ does not satisfy $(HS(F^{\ast })),$ Theorem 6 does not
apply, and indeed its conclusion fails. Also Theorem 2 (on linearity) --
which is at the heart of Theorem 6 (via Theorem 1) -- fails, as here the
domain of $G$ is $\mathbb{Q}$ and $G=0$ with a linear extension $S=0$ to all
of $\mathbb{R}.$ This shows the full force of Prop. 1(iv).

\bigskip

Variants of Th. 6 are possible. In the preceding argument Theorem 1 may be
replaced by Theorem 1$_{a}^{\prime }$ to yield:

\bigskip

\noindent \textbf{Theorem 6}$^{\prime }$\textbf{\ (Quantifier-Weakening
Theorem).} \textit{Theorem 6 holds with }(ii)\textit{\ replaced by:}

\noindent (ii\_a)$^{\prime }$\textit{\ }$F^{\ast }$\textit{\ satisfies the
Heiberg-Seneta boundedness condition thinned to a symmetric set }$\Sigma $%
\textit{\ that is locally SW, i.e.}%
\[
\lim \sup_{u\rightarrow 0,\text{ }u\in \Sigma }S(u)\leq 0. 
\]

Likewise, using Theorem 1$_{d}^{\prime }$ yields:

\bigskip

\noindent \textbf{Theorem 6}$^{\prime \prime }$\textbf{\ (Strong
Quantifier-Weakening Theorem). }\textit{Theorem 6 holds with }(ii)\textit{\
replaced by:}\newline
\noindent (ii\_d)$^{\prime }$\textit{\ }$F^{\ast }$\textit{\ is bounded on a
subset of }$\mathbb{A}_{F}$\textit{\ that is shift-compact (e.g. on a set
that is locally SW, and so on an open set).}

\bigskip

\noindent \textbf{Proof. }Again $\mathbb{A}_{F}$ is a subgroup and $%
G=F^{\ast }|\mathbb{A}_{F}$ is additive (Prop. 1); the condition $(\lim \sup
)_{0}$ follows from the assumption that $SW$-$\lim \sup {}_{M}(F^{\ast })$
holds for some $M.$ As $F^{\ast }$ is bounded on a shift-compact subset of $%
\mathbb{A}_{F},$ Theorem 1$^{\prime }$ applies to $S=F^{\ast }.$ $\square $

\bigskip

\noindent \textbf{Cautionary Example Again. }Recall that for $F:=\mathbf{1}_{%
\mathbb{R}\backslash \mathbb{Q}}$ one has $F^{\ast }=\mathbf{1}_{\mathbb{R}%
\backslash \mathbb{Q}}$ and $\mathbb{A}_{F}=\mathbb{Q}$. So here $\mathbb{A}%
_{F}$ is dense but not shift-compact, and $F^{\ast },$ though linear on $%
\mathbb{A}_{F},$ does not satisfy (ii).

\bigskip

The density assumption above may be weakened by using Theorem 4:

\bigskip

\noindent \textbf{Theorem 7.} \textit{With }$\mathbb{A}_{F}$ \textit{as
above,} \textit{suppose that there exists a shift-symmetric, analytic,
spanning set }$T\subseteq \mathbb{A}_{F}\ $\textit{such that }$F^{\ast }|T$%
\textit{\ is locally bounded above.}\newline
\noindent \textit{Then }$F^{\ast }$\textit{\ is linear: }$F^{\ast
}(u)=\lim_{x\rightarrow \infty }[F(u+x)-F(x)]=cu,$\textit{\ for some} $c\in 
\mathbb{R}$\textit{\ and all }$u\in \mathbb{R}.$

\bigskip

\noindent \textbf{Proof. }As $\mathbb{A}_{F}$ is a subgroup and $G=F^{\ast }|%
\mathbb{A}_{F}$ is additive (by Prop. 1), apply Theorem 4 to $S=F^{\ast }.$ $%
\square $

\bigskip

\noindent \textbf{Remark. }The sharpenings of Theorem 6 make use of \textit{%
relatives} of the Heiberg-Seneta condition. Their formulation draws on
shift-compactness (and so sequential) properties of\ various `test sets', $T$
say. The classical development relies on the classic Steinhaus property,
more properly: the Steinhaus-Weil\footnote{%
Here, as with Steinhaus, the context is $\mathbb{R}$; Weil's context is
(Haar) measurability in locally compact groups [Wei], cf. [GroE].} (interior
points) property of test sets (that $T-T$ has interior points): see BGT Th.
1.1.1; we study the links between the Steinhaus-Weil property and
shift-compactness elsewhere [BinO17].

\bigskip

The classical \textit{Quantifier Weakening Theorems} of regular variation
(BGT\ \S 1.4.3 and \S 3.2.5) are re-stated below as Theorems K and BKdH.
There, one needs as side-condition the Heiberg-Seneta condition $HS$
restated multiplicatively here as ($\lim \sup $) (or a thinned version of
it, as in Theorems 6$^{\prime }$, 6$^{\prime \prime }$). Recall from above
the $^{\ast }$ notation (as in $g^{\ast })$ signifying that limsup replaces $%
\lim .$

\bigskip

\noindent \textbf{Theorem K} (cf. BGT: Th. 1.4.3). \textit{Suppose that}%
\begin{equation}
\lim \sup\nolimits_{\lambda \downarrow 1}g^{\ast }(\lambda )\leq 1. 
\tag{$\lim \sup $}
\end{equation}%
\textit{Then the following are equivalent:}\newline
\textit{\noindent }(i)\textit{\ there exists }$\rho \in \mathbb{R}$\textit{\
such that }%
\[
f(\lambda x)/f(x)\rightarrow \lambda ^{\rho }\qquad (x\rightarrow \infty
)(\forall \lambda >0); 
\]%
\textit{\noindent }(ii)\textit{\ }$g(\lambda )=\lim_{x\rightarrow \infty
}f(\lambda x)/f(x)$\textit{\ exists, finite for all }$\lambda $\textit{\ in
a non-negligible set;}\newline
\textit{\noindent }(iii)\textit{\ }$g(\lambda )$\textit{\ exists, finite,
for all }$\lambda $\textit{\ in a dense subset of }$(0,\infty );$\newline
\textit{\noindent }(iv)\textit{\ }$g(\lambda )$\textit{\ exists, finite for }%
$\lambda =\lambda _{1},\lambda _{2}$\textit{\ with }$(\log \lambda
_{1})/\log \lambda _{2}$ \textit{irrational.}

\bigskip

Theorem K\ is an immediate corollary of Theorem 6, as (limsup) iff ($%
HS(F^{\ast })$).

\bigskip

\noindent \textbf{Theorem BKdH} (cf. BGT: Th. 3.2.5). \textit{For }$g$ 
\textit{with }%
\[
\lim_{x\rightarrow \infty }g(\lambda x)/g(x)=\lambda ^{\rho }\qquad (\lambda
>0), 
\]%
\textit{\ and}%
\begin{equation}
\lim \sup\nolimits_{\lambda \downarrow 1}f^{\ast }(\lambda )\leq 0, 
\tag{$\lim \sup $}
\end{equation}%
\textit{the following are equivalent:}\newline
\textit{\noindent }(i)\textit{\ }$k(\lambda ):=\lim_{x\rightarrow \infty }$[$%
f(\lambda x)-f(x)]/g(x)$\textit{\ exists, finite for all }$\lambda >0,$ 
\textit{and }$k(\lambda )=c[\lambda ^{\rho }-1]/\rho $\textit{\ for some }$c$%
\textit{\ and all }$\lambda $\textit{\ on a non-negligible set;}\newline
\textit{\noindent }(ii)\textit{\ }$k(\lambda )$ \textit{exists, finite for
all }$\lambda $\textit{\ in a non-negligible set;}\newline
\textit{\noindent }(iii)\textit{\ }$k(\lambda )$\textit{\ exists, finite,
for all }$\lambda $\textit{\ in a dense subset of }$(0,\infty );$\newline
\textit{\noindent }(iv)\textit{\ }$k(\lambda )$\textit{\ exists, finite for }%
$\lambda =\lambda _{1},\lambda _{2}$\textit{\ with }$(\log \lambda
_{1})/\log \lambda _{2}$ \textit{irrational.}

\bigskip

Implicit in the proofs in BGT is the \textit{Goldie functional equation }%
(GFE) and \textit{Goldie functional inequality }(GFI). This is made explicit
in [BinO14]. See (\textrm{Add}$_{\mathbb{A}}$); cf. BGT,\ Equation (3.2.7)),
and [Ost3] (on the relation between (GFE) and homomorphisms); we refer to
these sources for background. In results of this type, the usual
Baire/measurable assumptions are conspicuous by their absence. GFE\ in its
simplest form below bears little relation to CFE:%
\begin{eqnarray}
K(u+v) &=&K(u)+e^{u}K(v)\qquad (u,v\in \mathbb{A}),  \TCItag{$GFE$} \\
K(u+v) &\leq &K(u)+e^{u}K(v).  \TCItag{$GFI$}
\end{eqnarray}%
It is immediate from $(GFE)$ that either $K$ is trivial: $K\equiv 0,$ or for
some $\rho \neq 0$%
\[
K(u)=(e^{u}-1)/\rho \qquad (u\in \mathbb{A}). 
\]%
Only the latter case is of interest here. Define $H$ by%
\[
H(x):=K(\log x)\qquad (x\in \mathbb{E}:\mathbb{=}\exp \mathbb{A}). 
\]%
Writing $u=\log x$ etc. gives%
\begin{eqnarray*}
K(u+v) &=&K(\log (xy))=K(\log x)+xK(\log y), \\
H(xy) &=&H(x)+xH(y).
\end{eqnarray*}%
As $K(\log x)=(x-1)/\rho $ for some $\rho >0$, $H:$ $\mathbb{E\rightarrow R}$
is injective and order-preserving, so that $\mathbb{G}=K[\mathbb{A}]$ is
dense, if $\mathbb{A}$ is. Put 
\[
\eta (y):=1+\rho y, 
\]%
so that $H^{-1}(y)=\eta (y)$ for $y\in \mathbb{G}$. Now consider on $\mathbb{%
G}$ the Popa `circle' operation (Popa [Pop] in 1965, and Javor [Jav] in
1968):%
\[
x\circ y=x\circ _{\rho }y:=x+\eta (x)y. 
\]%
This is indeed a group operation with neutral element $1_{\rho }=0$ and
inverse $x_{\circ }^{-1}=-x/\eta (x);$ for background see [BinO15] and
[Ost3]. This group structure allows the Goldie equation to express
homomorphy: 
\[
H(xy)=H(x)\circ H(y)=H(x)+\eta (H(x))H(y)=H(x)+xH(y)\qquad (x,y\in \mathbb{E}%
). 
\]%
Alternatively, again as $H$ is invertible, with $X=H(x)$ etc, the equation
may be re-configured to the celebrated \textit{Go\l \k{a}b-Schinzel equation}%
:%
\begin{equation}
\eta (X\circ Y)=\eta (X)\eta (Y)\qquad (X,Y\in \mathbb{G}),  \tag{$GS$}
\end{equation}%
introduced for the study of one-parameter subgroups of affine groups, for
which see [AczD, Ch. 19] and the more recent [Brz].

Theorems 0, 0$^{\prime },$ and 0$^{+}$ are transferable, on the basis of two
simple facts stated in Prop. 12 below, to the context of functions $f:$ $%
\mathbb{R}_{+}\mathbb{\rightarrow G}$, with $\mathbb{G}$ equipped with the
Popa circle operation and with subadditivity replaced by:%
\[
f(xy)\leq f(x)\circ _{\rho }f(y) 
\]%
so as to yield GFI. (This development with its associated side-conditions
complements the alternative inequality 
\[
f(x+f(x)y)\leq f(x)f(y), 
\]%
studied in [Jab2]; cf. the `suboperative' functions of [HilP, 8.9].) Thus
the BKdH-RV version comes at little cost, as do the corresponding Quantifier
Weakening Theorems (Th. 6, 6$^{\prime }$, 6$^{\prime \prime }$) with linear $%
\rho x$ replaced by the affine $1+\rho x.$ We write 
\[
\mathbb{G}_{+}^{\rho }:=\{x:1+\rho x>0\}. 
\]

\bigskip

\noindent \textbf{Proposition 12. }\textit{For }$\rho \geq 0,$ \textit{the
set }$[0,\infty )$\textit{\ is a sub-semigroup of }$\mathbb{G}_{+}^{\rho };$%
\textit{\ the induced order: }$y\leq _{\rho }x$\textit{\ iff }$x\circ _{\rho
}y^{-1}\in \lbrack 0,\infty )$\textit{\ coincides with }$y\leq x.$ \textit{%
Furthermore, if }$c>0$\textit{\ and }$a<b,$ \textit{then}%
\[
a\circ _{\rho }c\leq b\circ _{\rho }c; 
\]%
\textit{in particular,}%
\[
(a,b)\circ _{\rho }c=(a\circ _{\rho }c,b\circ _{\rho }c), 
\]%
\textit{i.e. the Euclidean topology on }$\mathbb{R}_{+}$\textit{\ is
invariant under positive translation under }$\circ _{\rho }.$

\noindent \textit{Likewise, for }$\rho >0,$ \textit{if }$0<c<d,$\textit{\
and }$a<b$\textit{\ with }$a,b\in \mathbb{G}_{\rho }^{+},$ \textit{then}%
\[
a\circ _{\rho }c\leq b\circ _{\rho }d. 
\]

\bigskip

\noindent \textbf{Proof. }For the first assertion observe that 
\[
0\leq x-(1+\rho x)y/(1+\rho y)\text{ iff }0\leq x(1+\rho y)-(1+\rho x)y=x-y,%
\text{ as }1+\rho y>0. 
\]%
The rest is immediate, since, $\eta _{\rho }$ is positive on $\mathbb{G}%
_{+}^{\rho }\supseteq \mathbb{R}_{+}$ and order-preserving for $\rho >0$. $%
\square $

\bigskip

We hope to give a more detailed account of this area elsewhere. Suffice it
to say here that the two forms, Karamata and Bojani\'{c}-Karamata/de Haan,
of regular variation above are related to -- and indeed, subsumed within --
a third form, \textit{Beurling }slow and regular variation; see [BinO13,
10.3], [BinO15, \S 7]. The link here is the Popa circle operation above.

\section{ Complements}

\noindent 8.1 \textit{Sources. }This paper is a sequel to three sources:
[BinO11] (on the theorems of Steinhaus and Ostrowski, which underpin
everything in this area), [BinO14] (on the \textit{Goldie functional equation%
} (GFE) and inequality) and [Ost3] (on the relation between the more general 
\textit{Goldie-Beurling equation} and homomorphisms between Popa groups)
which here specializes to (GFE).

\noindent 8.2 \textit{Frullani integrals. }A classical situation where
quantifier weakening is important (which, as it happened, was the original
motivation for [BinG, I,II.6], and so for this paper and [BinO11] via BGT)
is the theory of \textit{Frullani integrals }([BinG, II \S 6], BGT \S 1.6.4;
Berndt [Bern])\textit{, }important in many areas of analysis and
probability. This in turn is a combination of two results from regular
variation, the Aljan\v{c}i\'{c}-Karamata theorem (a result of \textit{%
Mercerian} type) and the Characterisation Theorem (BGT \S 1.4), a central
result in the area inseparable from quantifier weakening. One reason why
regular variation is so ubiquitous and useful is its relevance to \textit{%
scaling }[Bin].

\noindent 8.3 \textit{Shift-compactness and Theorem S.} Evidently any
(non-degenerate) interval is shift-compact; more generally, so are
non-negligible Baire/measurable sets -- this is the \textit{%
Kestelman-Borwein-Ditor Theorem}, KBD, for which see [BinO11, Th. 4.2]. Any
shift-compact set $\Sigma $ has the \textit{classic Steinhaus property}
(terminology of [BarFN]): $0$ is an interior point of $\Sigma -\Sigma ,$ see
[BinO9, Th. 2]. The following combinatorial version of the Steinhaus
Subgroup Theorem, Theorem S\ below, will be seen capable of bearing the
burden of the proof of our version Theorem 1$_{d}^{\prime }$ above. See
[BinO11] for further equivalences in Theorem S (e.g. that $S$ has finite
index in $\mathbb{R},$ and statements involving Ramsey theory).

\bigskip

\noindent \textbf{Theorem S }([BinO11, Th. 6.2]) \textit{For an additive
subgroup }$\mathbb{A}$\textit{\ of }$\mathbb{R},$\textit{\ the following are
equivalent:}

\noindent (i) $\mathbb{A}=\mathbb{R},$

\noindent (ii) $\mathbb{A}$ \textit{contains a subset that is locally
Steinhaus-Weil (e.g. a non-negligible Baire/measurable set),}

\noindent (iii) $\mathbb{A}$ \textit{is shift-compact.}

\bigskip

\noindent 8.4 \textit{Bitopological Berz Theorem.} The proof of Theorem 3 in 
\S 5 above can be dualized to yield a parallel alternative and new proof for
Theorem B. Here, in place of Baire's continuity theorem, a careful use of
Lusin's theorem ([Hal, \S 55]; cf. [BinO8, \S 2] for a `near-analogue')
demonstrates linearity on a subset $\Sigma \cap \mathbb{R}_{+}$ covering
almost all of $\mathbb{R}_{+},$ and likewise on a subset $\Sigma \cap 
\mathbb{R}_{-}$ covering almost all of $\mathbb{R}_{-}$; then Props 10 and 7
above complete the proof. However, an argument proving simultaneously
Theorem 3 and Theorem B can be given [BinO16], by appeal to\textit{\
density-topology} arguments, for which see [BinO10,11], cf. [Wil] and [Ost1].

\noindent 8.5 \textit{Dependence on axioms of set theory. }For a summary of
the background information needed to appreciate the various set-theoretic
axioms which implicitly confront analysts we refer to Appendix 1 of the
fuller arXiv version of [BinO16]; the earlier article [Wri] of 1977 had a
similar motivation. This may be omitted by the expert (or uninterested)
reader.

\noindent 8.6 \textit{Kominek's Theorem. }We include this (discussed in \S %
6) here, as it is an immediate corollary of Prop. 6.

\bigskip

\noindent \textbf{Kominek's Theorem }([Kom], cf. [Jon]\textbf{). }\textit{%
For additive }$A:\mathbb{R}\rightarrow \mathbb{R}$, \textit{if }$A|T$\textit{%
\ is bounded on an analytic spanning set }$T$\textit{, then }$A$\textit{\ is
continuous.}

\bigskip

\noindent \textbf{Proof. }If $T$ is analytic and spans $\mathbb{R}$ then, as 
$A(T-T)$ is bounded, w.l.o.g. $T=-T$ (otherwise repace $T$ by $T\cup (-T)),$
and%
\[
\mathbb{R}=\bigcup\nolimits_{n\in \mathbb{N}}\bigcup\nolimits_{\mathbf{q}\in 
\mathbb{Q}^{n}}(q_{1}T+...+q_{n}T). 
\]%
So there are $n\in \mathbb{N}$ and $m_{1},...m_{n}\in \mathbb{N}$ with $%
S:=m_{1}T+...+m_{n}T$ of positive measure. So $A$ is bounded on $S+S$ and so
on an interval. Now apply Prop. 6 (with $\mathbb{A}=\mathbb{R}$). $\square $

\noindent 8.7 \textit{Kingman's Subadditive Ergodic Theorem. }Detailed study
of subadditivity is partially motivated by links with the Kingman
subadditive ergodic theorem, for which see e.g. [Kin1, 2], Steele [Ste].

\noindent 8.8 \textit{Sublinearity and risk measures. }An important class of
functions with the two properties of subadditivity and positive homogeneity
but with a more general domain occurs in mathematical finance -- the \textit{%
coherent risk measures} introduced by Artzner et al. [ArtDEH]; for textbook
treatments see [McNFE], [FolS, 4.1]. For the more general domains (and brief
commentary on the context) see again [BinO16].

\bigskip

\noindent \textbf{Acknowledgement. }We are most grateful to the Referee for
his careful reading and very thoughtful and challenging comments, which have
helped to improve this paper significantly, and to the Editor, Jan
Wiegerinck, for all his help.

\bigskip

\noindent \textbf{Postscript. }We take pleasure in dedicating this paper to
Charles Goldie. The first author is happy to recall that the argument in
[BinG] (and later in BGT) that gave rise to this paper was due to him, and
was the beginning of their long and fruitful collaboration.

\begin{center}
\textbf{References}
\end{center}

\noindent \lbrack AczD] J. Acz\'{e}l, J. Dhombres, \textsl{Functional
equations in several variables. With applications to mathematics,
information theory and to the natural and social sciences.} Encyclopedia of
Math. and its App. \textbf{31}, CUP, 1989\newline
\noindent \lbrack BarFN] A. Bartoszewicz, M. Filipczak, T. Natkaniec, On Sm%
\'{\i}tal properties. \textsl{Topology Appl.} \textbf{158} (2011),
2066--2075.\newline
\noindent \lbrack Bern] B. C. Berndt, Ramanujan's quarterly reports. \textsl{%
Bull. London Math. Soc.} \textbf{16} (1984), no. 5, 449--489.\newline
\noindent \lbrack Berz] E. Berz, Sublinear functions on $\mathbb{R},$ 
\textrm{\textsl{Aequat. Math.}} \textbf{12} (1975), 200-206.\newline
\noindent \lbrack Bin] N. H. Bingham, Scaling and regular variation. \textsl{%
Publ. Inst. Math. Beograd} \textbf{97} (111) (2015), 161-174.\newline
\noindent \lbrack BinG] N. H. Bingham, C.M. Goldie, Extensions of regular
variation. I. Uniformity and quantifiers. \textsl{Proc. London Math. Soc.}
(3) \textbf{44} (1982), 473--496; Extensions of regular variation. II.
Representations and indices. \textsl{Proc. London Math. Soc.} (3) \textbf{44}
(1982), 497--534.\newline
\noindent \lbrack BinGT] N. H. Bingham, C. M. Goldie and J. L. Teugels, 
\textsl{Regular variation}, 2nd ed., Cambridge University Press, 1989 (1st
ed. 1987). \newline
\noindent \lbrack BinO1] N. H. Bingham and A. J. Ostaszewski, Generic
subadditive functions, \textsl{Proc. Amer. Math. Soc. }\textbf{136} (2008),
4257-4266.\newline
\noindent \lbrack BinO2] N. H. Bingham and A. J. Ostaszewski, Beyond
Lebesgue and Baire: generic regular variation, \textsl{Colloquium Math.} 
\textbf{116} (2009), 119-138.\newline
\noindent \lbrack BinO3] N. H. Bingham and A. J. Ostaszewski, Infinite
combinatorics and the foundations of regular variation, \textsl{J. Math.
Anal. Appl.} \textbf{360} (2009), 518-529.\newline
\noindent \lbrack BinO4] N. H. Bingham and A. J. Ostaszewski, Automatic
continuity: subadditivity, convexity, uniformity\textrm{, \textsl{Aequat.
Math.}} \textbf{78} (2009), 257-270.\newline
\noindent \lbrack BinO5] N. H. Bingham and A. J. Ostaszewski, Infinite
combinatorics in function spaces: category methods, \textsl{Publ. Inst.
Math. Beograd} \textbf{86} (100) (2009), 55-73.\textrm{\ }\newline
\noindent \lbrack BinO6] N. H. Bingham and A. J. Ostaszewski, Automatic
continuity via analytic thinning\textrm{, \textsl{Proc. Amer. Math. Soc.}} 
\textbf{138} (2010), 907-919. \newline
\noindent \lbrack BinO7] N. H. Bingham and A. J. Ostaszewski, Normed groups:
Dichotomy and duality. \textsl{Dissert. Math.} \textbf{472} (2010), 138p. 
\newline
\noindent \lbrack BinO8] N. H. Bingham and A. J. Ostaszewski, Kingman,
category and combinatorics. \textsl{Probability and Mathematical Genetics}
(Sir John Kingman Festschrift, ed. N. H. Bingham and C. M. Goldie), 135-168,
London Math. Soc. Lecture Notes in Mathematics \textbf{378}, CUP, 2010. 
\newline
\noindent \lbrack BinO9] N. H. Bingham and A. J. Ostaszewski, Regular
variation without limits, \textsl{J. Math. Anal. Appl.} \textbf{370} (2010),
322-338.\newline
\noindent \lbrack BinO10] N. H. Bingham and A. J. Ostaszewski, Beyond
Lebesgue and Baire II: Bitopology and measure-category duality. \textsl{%
Colloquium Math.}, \textbf{121} (2010), 225-238.\newline
\noindent \lbrack BinO11] N. H. Bingham and A. J. Ostaszewski, Dichotomy and
infinite combinatorics: the theorems of Steinhaus and Ostrowski. \textsl{%
Math. Proc. Camb. Phil. Soc.} \textbf{150} (2011), 1-22. \newline
\noindent \lbrack BinO12] N. H. Bingham and A. J. Ostaszewski, Steinhaus
theory and regular variation: De Bruijn and after. \textsl{Indag. Math.} (N.
G. de Bruijn Memorial Issue), \textbf{24} (2013), 679-692.\newline
\noindent \lbrack BinO13] N. H. Bingham and A. J. Ostaszewski, Beurling slow
and regular variation, \textsl{Trans. London Math. Soc., }\textbf{1} (2014),
29-56 (fuller versions: arXiv 1301.5894 and 1307.5305).\newline
\noindent \lbrack BinO14] N. H. Bingham and A. J. Ostaszewski, Cauchy's
functional equation and extensions: Goldie's equation and inequality, the Go%
\l \k{a}b-Schinzel equation and Beurling's equation, \textsl{Aequationes
Math.}, \textbf{89.5} (2015), 1293-1310, arXiv1405.3947.\newline
\noindent \lbrack BinO15] N. H. Bingham and A. J. Ostaszewski, Beurling
moving averages and approximate homomorphisms, \textsl{Indag. Math. }\textbf{%
27} (2016), 601-633 (fuller version: arXiv1407.4093).\newline
\noindent \lbrack BinO16] N. H. Bingham and A. J. Ostaszewski,
Category-measure duality: convexity, mid-point convexity and Berz
sublinearity, \textsl{Aequationes Math.}, \textbf{91.5} (2017), 801--836, (
fuller version: arXiv1607.05750).\newline
\noindent \lbrack BinO17] N. H. Bingham and A. J. Ostaszewski, Beyond
Lebesgue and Baire IV: Density topologies and a converse Steinhaus-Weil
theorem, \textsl{Topology and its Applications}, to appear, arXiv1607.00031.%
\newline
\noindent \lbrack BinO18] N. H. Bingham and A. J. Ostaszewski, The
Steinhaus-Weil property: its converse, Solecki amenability and
subcontinuity, arXiv1607.00049.\newline
\noindent \lbrack BojK] R. Bojani\'{c} and J. Karamata, \textsl{On a class
of functions of regular asymptotic behavior, }Math. Research Center Tech.
Report 436, Madison, Wis. 1963; reprinted in \textsl{Selected papers of
Jovan Karamata} (ed. V. Mari\'{c}, Zevod za Ud\v{z}benike, Beograd, 2009),
545-569.\newline
\noindent \lbrack Brz] J. Brzd\k{e}k, The Go\l \k{a}b-Schinzel equation and
its generalizations, \textsl{Aequat. Math.} \textbf{70} (2005), 14-24.%
\newline
\noindent \lbrack CrnGH] M. Crnjac, B. Gulja\v{s}, H. I. Miller, On some
questions of Ger, Grubb and Kraljevi\'{c}. \textsl{Acta Math. Hungar.} 
\textbf{57} (1991), 253--257.\newline
\noindent \lbrack Darb] G. Darboux, Sur la composition des forces en
statique, \textsl{Bull. Sci. Math.} \textbf{1} (9) (1875), 281-288.\newline
\noindent \lbrack Darj] U. B. Darji, On Haar meager sets. \textsl{Topology
Appl.} \textbf{160} (2013), 2396--2400.\newline
\noindent \lbrack Ful] R. V. Fuller, Relations among continuous and various
non-continuous functions. \textsl{Pacific J. Math.} \textbf{25} (1968),
495--509.\newline
\noindent \lbrack GroE] K.-G. Grosse-Erdmann, An extension of the
Steinhaus-Weil theorem. \textsl{Colloq. Math.} \textbf{57}.2 (1989),
307--317.\newline
\noindent \lbrack dHa] L. de Haan, \textsl{On regular variation and its
applications to the weak convergence of sample extremes.} Math. Centre
Tracts 32, Amsterdam 1970.\newline
\noindent \lbrack Hal] P. R. Halmos, \textsl{Measure Theory}, Grad. Texts in
Math. 18, Springer 1974. (1st. ed. Van Nostrand, 1950).\newline
\noindent \lbrack HarW] G. H. Hardy and E. M. Wright, An introduction to the
theory of numbers, 6th ed. Revised by D. R. Heath-Brown and J. H. Silverman
(OUP 2008).\newline
\noindent \lbrack Hei] C. H. Heiberg, A proof of a conjecture by Karamata. 
\textsl{Publ. Inst. Math. (Beograd)} (N.S.) \textbf{12} (26) (1971), 41--44.%
\newline
\noindent \lbrack HilP] E. Hille and R. S. Phillips, \textsl{Functional
analysis and semi-groups}, Coll. Publ. Vol 31, Amer. Math. Soc, 3rd ed. 1974
(1st ed. 1957).\newline
\noindent \lbrack Jab1] E. Jab\l o\'{n}ska, Some analogies between Haar
meager sets and Haar null sets in abelian Polish groups. \textsl{J. Math.
Anal. Appl.} \textbf{421} (2015),1479--1486.\newline
\noindent \lbrack Jab2] E. Jab\l o\'{n}ska, On solutions of a composite type
functional inequality,\textsl{\textsl{\ Math. Ineq. Apps. }}\textbf{18}
(2015), 207-215.\textsl{\textsl{\newline
}}\noindent \lbrack Jab3] E. Jab\l o\'{n}ska, A theorem of Piccard's type
and its applications to polynomial functions and convex functions of higher
orders. \textsl{Topology Appl.} \textbf{209} (2016), 46--55.\newline
\noindent \lbrack Jab4] E. Jab\l o\'{n}ska, A theorem of Piccard's type in
abelian Polish groups.\textsl{\ Anal. Math.} \textbf{42} (2016), 159--164.%
\newline
\noindent \lbrack Jav] P. Javor, On the general solution of the functional
equation $f(x+yf(x))=f(x)f(y).$ \textsl{Aequat. Math.} \textbf{1} (1968),
235-238.\newline
\noindent \lbrack Jon] F. B. Jones, Measures and other properties of Hamel
bases, \textsl{Bull. Amer. Math. Soc.} \textbf{48} (1942), 472-481.\newline
\noindent \lbrack Kin1] J. F. C. Kingman, Subadditive ergodic theory. 
\textsl{Ann. Probability} \textbf{1} (1973), 883--909.\newline
\noindent \lbrack Kin2] J. F. C. Kingman, Subadditive processes. \textsl{%
\'{E}cole d'\'{e}t\'{e} de Probabilit\'{e}s de Saint Fleur} \textsl{V}, 
\textsl{Lecture Notes in Math.} \textbf{539} (1976), 167-223.\newline
\noindent \lbrack Kom] Z. Kominek, On the continuity of $\mathbb{Q}$-convex
and additive functions, \textrm{\textsl{Aequat. Math.}} \textbf{23} (1981),
146-150.\newline
\noindent \lbrack Kuc] M. Kuczma, \textsl{An introduction to the theory of
functional equations and inequalities. Cauchy's equation and Jensen's
inequality.} 2nd ed., Birkh\"{a}user, 2009 [1st ed. PWN, Warszawa, 1985].%
\newline
\noindent \lbrack LukMZ] J. Luke\v{s}, J. Mal\'{y}, L. Zaj\'{\i}\v{c}ek, 
\textsl{Fine topology methods in real analysis and potential theory},
Lecture Notes in Mathematics \textbf{1189}, Springer, 1986.\newline
\noindent \lbrack MatS] J. Matkowski, T. \'{S}wi\k{a}tkowski, On subadditive
functions. \textsl{Proc. Amer. Math. Soc.} \textbf{119} (1993), 187--197.%
\newline
\noindent \lbrack MatZ] E. Mato\u{u}skov\'{a}, M. Zelen\'{y}, A note on
intersections of non--Haar null sets, \textsl{Colloq. Math.} \textbf{96}
(2003), 1-4.\newline
\noindent \lbrack Meh] M. R. Mehdi, On convex functions,\textrm{\ \textsl{J.
London Math. Soc.} }\textbf{39} (1964), 321-328.\newline
\noindent \lbrack MilO] H. I. Miller and A. J. Ostaszewski, Group actions
and shift-compactness. \textsl{J. Math. Anal. Appl.} \textbf{392} (2012),
23-39.\newline
\noindent \lbrack Ost1] A. J. Ostaszewski, Analytically heavy spaces:
Analytic Cantor and Analytic Baire Theorems, \textsl{Topology and its
Applications}, \textbf{158} (2011), 253-275.\newline
\noindent \lbrack Ost2] A. J. Ostaszewski, Beyond Lebesgue and Baire III:
Steinhaus' Theorem and its descendants, \textsl{Topology and its App.} 
\textbf{160} (2013), 1144-1154.\newline
\noindent \lbrack Ost3] A. J. Ostaszewski, Homomorphisms from Functional
Equations: The Goldie Equation, \textsl{Aequationes Math. }\textbf{90}
(2016), 427-448, (arXiv: 1407.4089).\newline
\noindent \lbrack Ostr] A. Ostrowski\textrm{\textsc{,} }Mathematische
Miszellen XIV: \"{U}ber die Funktionalgleichung der Exponentialfunktion und
verwandte Funktionalgleichungen,\textrm{\ \textsl{Jahresb. Deutsch. Math.
Ver.} }\textbf{38} (1929) 54-62, reprinted in\textrm{\ \textsl{Collected
papers of Alexander Ostrowski}, }Vol. 4, 49-57, Birkh\"{a}user, Basel, 1984. 
\newline
\noindent \lbrack Oxt] J. C. Oxtoby: \textsl{Measure and category}, 2nd ed.
Graduate Texts in Math. \textbf{2}, Springer, 1980.\newline
\noindent \lbrack Pop] C. G. Popa, Sur l'\'{e}quation fonctionelle $%
f[x+yf(x)]=f(x)f(y),$ \textsl{Ann. Polon. Math.} \textbf{17} (1965), 193-198.%
\newline
\noindent \lbrack Roc] R. T. Rockafellar, \textsl{Convex analysis.}
Princeton, NJ, 1970 (reprinted 1997).\newline
\noindent \lbrack Rog] C. A. Rogers, J. Jayne, C. Dellacherie, F. Tops\o e,
J. Hoffmann-J\o rgensen, D. A. Martin, A. S. Kechris, A. H. Stone, \textsl{%
Analytic sets.} Academic Press, 1980.\newline
\noindent \lbrack Sen] E. Seneta, An interpretation of some aspects of
Karamata's theory of regular variation. \textsl{Publ. Inst. Math. (Beograd)}
(N.S.) \textbf{15 }(29) (1973), 111--119.\newline
\noindent \lbrack Sol] S. Solecki, Amenability, free subgroups, and Haar
null sets in non-locally compact groups. \textsl{Proc. London Math. Soc.} 
\textbf{(3) 93} (2006), 693--722.\newline
\noindent \lbrack Ste] J. M. Steele, Kingman's subadditive ergodic theorem, 
\textsl{Ann. de l'Institut Henri Poincar\'{e}}, \textbf{25} (1989), 93-98.%
\newline
\noindent \lbrack Wei] A. Weil, \textsl{L'integration dans les groupes
topologiques}, Actualit\'{e}s Scientifiques et Industrielles 1145, Hermann,
1965 (1$^{\text{st }}$ ed. 1940).\newline
\noindent \lbrack Wil] W. Wilczy\'{n}ski, \textsl{Density topologies},
Handbook of measure theory, Vol. I, II, 675--702, North-Holland, 2002.%
\newline
\noindent \lbrack Wri] J. D. Maitland Wright, Functional Analysis for the
practical man, 283--290 in Functional Analysis: Surveys and Recent Results 
\textbf{27}, North-Holland Math. Studies, 1977.

\bigskip

\bigskip

\noindent Mathematics Department, Imperial College, London SW7 2AZ;
n.bingham@ic.ac.uk \newline
Mathematics Department, London School of Economics, Houghton Street, London
WC2A 2AE; A.J.Ostaszewski@lse.ac.uk\newpage

\end{document}